\documentclass[12pt,reqno]{amsart}
\usepackage{mathrsfs}
\usepackage{amssymb}
\usepackage{amsfonts}
\usepackage{amsbsy}
\usepackage{latexsym}
\usepackage{amssymb,latexsym,amsmath,amsthm}
\usepackage{framed}
\usepackage{subcaption}
\usepackage{booktabs}
\usepackage{diagbox}
\usepackage{enumitem}

\DeclareCaptionLabelFormat{tightfig}{#1\hspace{0.15em}#2}
\usepackage[colorlinks,linkcolor=blue]{hyperref}
\usepackage{graphicx}
\usepackage{float}
\usepackage{xcolor}
\usepackage{epstopdf}
\usepackage{bm}
\usepackage{cite}
\usepackage{makecell}
\usepackage{multirow}
\usepackage{placeins}
\let\mathbf\bm
\DeclareFontFamily{U}{mathx}{}
\DeclareFontShape{U}{mathx}{m}{n}{<-> mathx10}{}
\DeclareSymbolFont{mathx}{U}{mathx}{m}{n}
\DeclareMathAccent{\widecheck}{0}{mathx}{"71}
\theoremstyle{plain}

 \theoremstyle{remark} 

\newcommand{\Z}{{\mathbb Z}}

\newcommand{\N}{{\mathbb N}}

\newtheorem {theo} {\bf Theorem} [section]
\newtheorem {prop} [theo] {\bf Proposition}

\newtheorem {lem} [theo] {\bf Lemma}

\newtheorem {defi} {\bf Definition}[section]
\newtheorem{exam} {\bf Example}[section]

\newtheorem{rem}{\bf Remark}[section]

\numberwithin{equation}{section}
\begin{document}
	\title[Characterization of AV-LT and  sample-based reconstruction]{Characterization of  Attenuated V-Line Transform and Direct Reconstruction from Samples}
	\author{Youfa Li}
	\address{School of Mathematics \& Center for Applied Mathematics of Guangxi (Guangxi University) \\
		Guangxi University,  Nanning, China }
	\email{youfalee@hotmail.com}
	\author{Lvchang Li}
	\address{School  of Mathematics\\
		Guangxi University,  Nanning, China }
	\email{llc0118@126.com}
	\date{\today}

	\begin{abstract} 
		The V‑line transform (V‑LT), a generalized Radon transform (RT),  is widely used  in Compton scattering  imaging (CSI).
		The central  mathematical problem  of CSI is to reconstruct a function from its V-LT.
		Unlike V-LT, the attenuated V‑line transform (AV-LT)  accounts for the  medium‑induced  loss of intensity, and has been demonstrated to improve image quality.
		The present  paper makes   two main contributions. \textbf{(1)} We establish the fundamental  characterizations for AV-LT including  operator properties ($L^{2}$-boundedness, operator norm etc.), a necessary and sufficient condition  for the preservation of compact support,  and sharp estimation of  Sobolev smoothness. These   
		characterizations remain   unexplored, and they  reveal   essential differences among    AV-LT,  RT and V-LT.
		\textbf{(2)} We develop a direct  reconstruction method based on   AV-LT samples that does not require solving any   equation. Our sampling formula gives exact reconstruction for band-limited functions. We also provide  approximate reconstruction for   functions in $H^s(\mathbb R^2)$ with $s>2$. Our reconstruction method is \emph{different} from the existing ones.  Firstly, the existing methods fall mainly  into two categories. They either use inversion formulas requiring the full AV-LT data or solve often ill-posed inverse problems from samples via regularization. Instead, our 
		method is a direct reconstruction that is  based on AV-LT samples  and does not need to solve any equation.
		Secondly, the existing  methods commonly  assume infinitely differentiable functions. In contrast, our   applicable 
		function space $H^s(\mathbb R^2)$ contains many 
		functions that are not infinitely differentiable. 
		Numerical simulations are conducted to verify the results. Particularly, the simulations show that our direct reconstruction method outperforms the regularization method in terms of robustness to noise, when reconstructing band-limited functions.

	\end{abstract}
	\maketitle
	\section{Introduction}
	\subsection{Background}
	Relying on the detection of gamma-ray photons,  single-photon emission computed tomography (SPECT) is an advanced medical  diagnosis technique for  functional imaging. As a promising substitute for the Anger camera for  gamma‑ray imaging, the Compton camera can greatly improve the sensitivity of SPECT (cf. \cite{HaltmeierMoonSchiefeneder}). The Compton camera-based imaging (reconstruction of functions) is achieved by the inversion of V-line transform (V-LT, a generalized Radon transform, cf. \cite{Ambartsoumian—shu}). However, the  attenuation of photons is neglected in the  reconstruction method based on V-LT, which may cause notable degradation of image quality.  Haltmeier,  
	Moon  and  Schiefeneder \cite{HaltmeierMoonSchiefeneder} introduced 
	the attenuated V‑line transform (AV-LT) and used it to improve the image quality. In this paper, we are interested in the mathematical characterizations of  AV-LT and the reconstruction problem based on AV-LT samples. Before presenting our motivations, we first recall the definitions  of V-LT and AV-LT.


	\subsubsection{V-line transform and attenuated V‑line transform}    
\begin{figure}[htbp]
		\centering
		\includegraphics[width=0.38\textwidth]{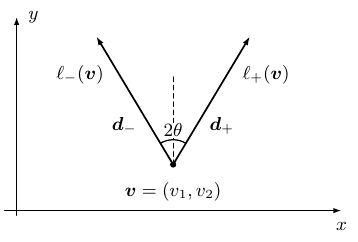}
		\caption{The V-line/ AV-line geometry.}
		\label{fig:vline-geometry}
	\end{figure}

	Associated  with the  half-opening angle  $\theta\in(0,\pi/2)$, the V-LT (cf. \cite{Ambartsoumian2016})  of  $f\in L^2(\mathbb{R}^2)$  at  $\bm{x}\in \mathbb{R}^2$  is defined to be  the sum of  the integrals along the two rays  $\ell_s(\mathbf{x}):=\{\mathbf{x}+t\mathbf{d}_s(\theta):t\ge0\}$ as follows
	\begin{equation}\label{Vline}
		(V_{\theta}f)(\mathbf{x})
		:=
		\sum_{s\in\{+,-\}}
		\int_0^\infty
		f\bigl(\mathbf{x}+t\mathbf{d}_s(\theta)\bigr)dt,
	\end{equation}
	where $\mathbf{d}_s(\theta):=(s\sin\theta,\cos\theta)\in\mathbb{S}^1$. Here $\bm{x}=(x_{1}, x_{2})$ is called the emanating point. 
	See Fig.\ref{fig:vline-geometry} for the geometric model for  V-LT. 

	
	By \cite{HaltmeierMoonSchiefeneder}, 
	the AV-LT is defined as 
	\begin{equation}\label{eq:Def-V-line-transform}
		(V_{\theta,\mu}f)(\mathbf{x})
		=\sum_{s\in\{+,-\}}\int_0^\infty f\bigl(\mathbf{x}+t\mathbf{d}_s(\theta)\bigr)e^{-\mu t}dt,
	\end{equation}
	with $\mu>0$ denoting the attenuation coefficient. For convenience  of narration,  we write $V_{\theta,\mu}f=T_+f+T_-f$,
	where
	\begin{equation}\label{eq:attenuated-ray-transforms-right}
		(T_+f)(\mathbf{x}):=\int_0^\infty f\bigl(\mathbf{x}+t\mathbf{d}_+(\theta)\bigr)e^{-\mu t}dt
	\end{equation}
	and
	\begin{equation}\label{eq:attenuated-ray-transforms-left}
		(T_-f)(\mathbf{x}):=\int_0^\infty f\bigl(\mathbf{x}+t\mathbf{d}_-(\theta)\bigr)e^{-\mu t}dt.
	\end{equation}
	See also  Fig.\ref{fig:vline-geometry} for 
	the geometric model for  AV-LT.

	\begin{rem}
		(1) Compared with V-LT, an attenuation term  $e^{-\mu t}$
		is introduced in AV-LT. (2) As two kinds of generalized Radon transforms, V-LT and   AV-LT are obtained by integrating the function $f$ along a V‑line. This is different from the Radon transform
		defined by 
		\begin{align}\label{Radonbianhuandedingyi}
			R_{\vartheta}f(\omega)
			:=
			\int_{x\cdot \vartheta=\omega} f(x)dx,
			\ 
			\omega\in\mathbb{R},\ \vartheta\in \mathbb{S}^1
		\end{align}
		which integrates along a straight line.
	\end{rem}
	
	Fixing the  half-opening angle  $\theta$  helps to simplify the scanning process (cf. \cite{Cebeiro2025}). Recently,  the inversion problems  for both  V-LT and   AV-LT with fixed  $\theta$ attract  considerable interest (e.g.\cite{GouiaZarradAmbartsoumian,AmbartsoumianLatifiJebelliMishra,AmbartsoumianMishraZamindar2024,Cebeiro2025,GouiaZarradMoon,JeonMoon,AmbartsoumianJebelli2019}).   Throughout the paper, we are interested in the model of fixed $\theta$. 
	
	
	
	
	\subsection{Motivations}



	\subsubsection{Properties  of  AV-LT are not trivial and  remain  unexplored} \label{somemotivations}
	Compared with V-LT and  RT, AV-LT exhibits some  essential differences in  the  operator  properties, preservation of  compact support  and  Sobolev smoothness.

	\begin{itemize}
		\item \underline{Operator properties}.
		We denote  the space of bounded linear operators on $L^2(\mathbb R^2)$
		by  $\mathcal L(L^2(\mathbb R^2),L^2(\mathbb R^2))$.
		It is well-known that RT does not lie in 
		$\mathcal L(L^2(\mathbb R^2),$
		$L^2(\mathbb R^2))$ (cf. \cite{Natterer2001}).
		For the convergence of $V_{\theta}f$, $f$ is required to satisfy additional conditions, including \textbf{sufficiently fast decay}  at infinity (cf. \cite{Ambartsoumian—shu}).
		This implies that $V_{\theta}$ may not  map $L^2(\mathbb{R}^2)$ into itself.
		In contrast to RT and  V-LT,  AV-LT introduces  an exponential decay factor $e^{-\mu t}$.
		One   naturally  asks the following   problems:  
		
		$\textbf{Q}1$:   Does  $V_{\theta, \mu}$ belong to $ \mathcal{L}(L^2(\mathbb{R}^2), L^2(\mathbb{R}^2))$? 
		
		$\textbf{Q}2$: If $V_{\theta, \mu}\in \mathcal{L}(L^2(\mathbb{R}^2), L^2(\mathbb{R}^2))$, what is the operator norm $\|V_{\theta, \mu}\|_{L^2(\mathbb R^2)\to L^2(\mathbb R^2)}$?
		

		\item \underline{Sobolev smoothness}.
		The Sobolev smoothness of $g\in L^2(\mathbb{R}^d)$ is defined as 
		\begin{align}\label{Sobolevguanghuadingyi}
			\nu_{2}(g)=\sup\{s\in \mathbb{R}: g\in H^{s}(\mathbb{R}^d)\}
		\end{align}
		where $H^{s}(\mathbb{R}^d)$ is the Sobolev space.  One of our purposes is to establish 
		the reconstruction  formula for $f\in L^{2}(\mathbb{R}^{2})$  from  the samples of   $V_{\theta, \mu}f$. 
		To make sampling  sense, $V_{\theta, \mu}f$  is required to be continuous   (cf. \cite{Aldroubi,Nashed,Sun1}).
		On the other hand, Sobolev smoothness is a typical criterion for continuity since $g\in L^2(\mathbb{R}^2)$ is continuous if  $\nu_{2}(g)>1$ (cf. \cite{LiHanSCICHINA}). Then one naturally asks:
		
		$\textbf{Q}3$:  What condition on $f\in L^{2}(\mathbb{R}^{2})$ ensures that $\nu_{2}(V_{\theta,\mu}f)>1$?

		For any  compactly supported $f\in L^2(\mathbb{R}^2)$,  it follows from  Han \cite[Lemma 2.4]{Han} 
		or L. H{\"o}rmander \cite[Theorem~4.3.1]{Hormander1971} 
		that 
		$\nu_{2}(R_{\vartheta} f) \geq \nu_{2}(f) $. Motivated by this, 
		we will investigate how AV-LT affect 
		the Sobolev smoothness.
		Based on  this, 
		we will give an answer to  $\textbf{Q}3$.

		\item \underline{Compact support}. 
		By Li, Fan and Han  \cite{LiFanHan}, $R_{\vartheta}f$ is compactly supported for any compactly supported $f\in L^{2}(\mathbb{R}^2)$. Then one  naturally asks:  
		
		$\textbf{Q}4$:  Is $V_{\theta, \mu}f$ also  compactly supported whenever  $f$ is compactly supported?  If not, what are the necessary and sufficient conditions on $f$ for $V_{\theta,\mu}f$ to be compactly supported?

	\end{itemize}	
	
	\subsubsection{Direct reconstruction formula from AV-LT samples  remains unexplored.}\label{xianyoufangfabingtai}
	The central  mathematical problem  for  Compton scattering  imaging (based on AV-LT) is to reconstruct a function $f$ from its AV-LT.  Current research on such a  reconstruction problem mainly  falls into two lines   as follows.  
	
	\textbf{(1) Inversion formula from AV-LT. }  
	Recently,  extensive research interest has been devoted  to establishing  inversion 
	formulas from AV-LT (e.g., \cite{GouiaZarradMoon
		,JeonMoon,HaltmeierMoonSchiefeneder,HaltmeierSchiefeneder2018,Duy}).
	Particularly, 
	Jeon and Moon \cite{JeonMoon} obtained a singular value decomposition (SVD) for an attenuated conical transform (with fixed half-opening angle),  and derived the corresponding inversion formula  in the  suitable Hilbert space. These inversion formulas  
	commonly   depend  on an integral involved  with AV-LT.  
	The    integral    requires all the values of AV-LT. 
	
	\textbf{(2) Reconstruction  from  AV-LT samples.} 
	In practice, what we   obtain are the samples of AV-LT but not all its values. One can 
	first give the representation  of the object $f$
	via  the prescribed    bases such as Fourier bases (e.g. \cite{HaltmeierMoonSchiefeneder}) and spherical harmonics (e.g.\cite{HaltmeierSchiefeneder2018}). Based on the representation, one  can reconstruct $f$ 
	(or, equivalently, the representation coefficients) from its AV-LT samples. 
	However, this reconstruction typically requires solving equations
	and suffers from severe ill-posedness
	(cf. \cite{HaltmeierMoonSchiefeneder,HaltmeierSchiefeneder2018}).
	
	Combining (1) and (2), one  naturally asks:
	
	$\textbf{Q}5$: Is it possible to develop a reconstruction formula, based on samples of the AV‑LT with a fixed half‑opening angle, that relies solely on evaluating the formula rather than solving any equation?
	
	This problem remains unexplored. From the perspective of sampling theory, we address $\textbf{Q}5$ for band-limited functions and for functions in $H^{s}(\mathbb{R}^{2})$, $s>2$, obtaining the  exact and approximate reconstruction, respectively.

	\subsection{Main contributions}\label{zhuyaogongxian123}
	
	The main contributions   of this paper are  organized as follows. 
	\begin{enumerate}

		\item\textbf{Operator properties: answers to  $\textbf{Q1-Q2}$.}
		For any $\mu>0$, we prove that $V_{\theta, \mu} \in \mathcal{L}(L^2(\mathbb{R}^2),$
		$ L^2(\mathbb{R}^2))$ and give its operator norm (Theorem \ref{L2-boundedness}). 

		\item\textbf{Sobolev smoothness of AV-LT and an answer to \textbf{Q3}.}
		(i) We establish  the sharp  estimate  $\nu_{2}(V_{\theta, \mu}f)\geq \nu_{2}(f)$ for any $f\in L^{2}(\mathbb{R}^{2})$. Here the  sharpness means that, there does not exist $\epsilon>0$ such that 
		$\nu_{2}(V_{\theta, \mu}f)\geq\nu_{2}(f) +\epsilon$
		for any $f\in L^{2}(\mathbb{R}^{2})$. 
		(ii) $\nu_{2}(f)>1$ ensures  $\nu_{2}(V_{\theta,\mu}f)>1$ in \textbf{Q3}. 
		All these results are organized in Theorem \ref{Vline-sobolev}. 
		
		
		\item\textbf{An answer to \textbf{Q4}.}
		(i) From the perspective of Fourier-Laplace transform,  we establish a necessary and sufficient condition for any compactly supported  $f\in L^{2}(\mathbb{R}^2)$, under which    $V_{\theta,\mu}f$ is also   compactly supported (Theorem \ref{thm:compact-support-V}).
		(ii) The   condition is not trivial since 
		$V_{\theta,\mu}f$ is not compactly supported for any  nonnegative and compactly supported   $0\not \equiv f\in L^2(\mathbb{R}^{2})$ (Theorem  \ref{compact-remark}).

		\item\textbf{An answer to \textbf{Q5} for   the band-limited function space.}
		We establish an exact reconstruction formula  for band-limited functions from the samples of their AV-LTs (Theorem \ref{thm:sampling-reconstruction}). 
		The key ingredient for establishing  the formula  is  to 
		\textcolor{black}{\textbf{design   the reconstruction  kernel}}. 
		The reconstruction   possesses two advantages: (i) it relies only  on samples of the AV‑LT instead of its full set of values; (ii) the   implementation requires only the evaluation of the formula rather than solving any equation.  
		

		
		\item \textbf{An answer to \textbf{Q5} for   Sobolev space.}
		We  apply  the  reconstruction scheme in (4)  to Sobolev spaces. Specifically,  for every $f\in H^s(\mathbb{R}^2)$ with $s>2$, we construct an  approximation $f_N^\#$ from samples of $V_{\theta,\mu}f$,  and establish an $L^2$-error estimate  which shows  that $f_N^\#$ converges exponentially   to $f$ as the level  $N\to\infty$ (Theorem \ref{thm:sobolev-rate}). Like the reconstruction formula in (4), $f_N^\#$
		also possesses the above  two advantages (i)
		and (ii). 
	\end{enumerate}
	
	\subsection{Highlights}
	\begin{enumerate}
		
		\item  \textbf{A different  reconstruction perspective from the existing methods.} 
		\begin{enumerate}
			\item 
			As mentioned in section \ref{xianyoufangfabingtai} (2),  the existing reconstruction methods   from AV-LT samples consists of two main steps: (i)  representing $f$ with respect to  the prescribed 
			bases; (ii)  reconstructing   $f$  reduces to reconstructing the   representation  coefficients. Step (ii) is an inverse problem and  commonly suffers from severe  ill‑posedness.

			\item 
			Unlike the existing methods, for the band-limited function $f$,  our reconstruction consists of two main steps: (i') By Shannon sampling theorem, 
			$V_{\theta,\mu}f$ is represented  by its samples; (ii') Based on this representation, we first  design the reconstruction kernel and then establish  the reconstruction formula. 
		\end{enumerate}
		
		\item \textbf{The advantages of the reconstruction formula.}  Our  reconstruction formulas mentioned  in section \ref{zhuyaogongxian123} (4, 5)  enjoy  three  advantages.  
		\begin{enumerate}
			\item They  require the samples but not all the values of AV-LT. This is different from the existing  inversion formulas. 
			\item They requires only evaluating the reconstruction formula  and  do  not require solving any equation. 
			This is different from the existing regularization-based methods. Simulations reveal that our direct reconstruction method outperforms the regularization approach in noise robustness, when recovering band-limited functions.
			\item    The existing AV-LT reconstruction methods require the object $f$ to be infinitely differentiable (e.g. \cite{Duy,GouiaZarradMoon,HaltmeierMoonSchiefeneder}). The reconstruction formula  mentioned in section \ref{zhuyaogongxian123} (5)
			holds for $H^s(\mathbb{R}^2)$ where it is only required that  $s>2.$ 
		\end{enumerate}

		\item   \textbf{We clarify the similarities and differences among RT, V-LT, and AV-LT}  in terms of $L^2$-boundedness, compact support preservation, and variation  in Sobolev smoothness.  They are summarized  in Table \ref{tab:RT-AVLT-comparison}.
	\end{enumerate}
	
	\begin{table}[!htbp]
		\centering
		\caption{Comparison of RT, AV-LT and V-LT w.r.t  $L^2$-boundedness, smoothness  and compact support.}
		\label{tab:RT-AVLT-comparison}
		\resizebox{\textwidth}{!}{%
			\begin{tabular}{lccc}
				\hline
				& \multirow{2}{*}{$L^2$-bounded?}
				& \multicolumn{2}{c}{Any compactly supported $f\in L^{2}(\mathbb{R}^{2})$} \\
				\cmidrule(lr){3-4}
				&
				& \shortstack[c]{Sobolev smoothness}
				& \shortstack[c]{compactly supported}? \\
				\hline
				RT
				& No (cf.\cite{Natterer2001})
				& \makecell[c]{$\nu_{2}(R_{\vartheta}f)\geq\nu_{2}(f)$ \ (cf.\cite{Han})\\$\nu_{2}(R_{\vartheta}f)>\nu_{2}(f)?$\\Not necessarily\\(Example \ref{example3.2})}
				& Yes \ (cf.\cite{LiFanHan})\\
				\hline
				AV-LT
				& \makecell[c]{Yes\\
					(Theorem \ref{L2-boundedness})}
				& \makecell[c]{$\nu_{2}(V_{\theta, \mu}f)\geq\nu_{2}(f)$\\(Theorem \ref{Vline-sobolev} (i))\\$\nu_{2}(V_{\theta, \mu}f)>\nu_{2}(f)?$\\Not necessarily\\(Theorem \ref{Vline-sobolev} (ii), Example \ref{ex:strict-sobolev})}
				& \makecell[c]{Not necessarily\\(Theorem \ref{thm:compact-support-V})}
				\\
				\hline
				V-LT
				& \makecell[c]{$\mathord{*}\mathord{*}\mathord{*}^{\textsuperscript{a}}$\\(Proposition \ref{Vlinebushoulian})}
				& \makecell[c]{$\mathord{*}\mathord{*}\mathord{*}^{\textsuperscript{a}}$\\(Proposition \ref{Vlinebushoulian})}
				& \makecell[c]{Not necessarily\\(Remark \ref{bumanzujinzhcheng})}
				\\
				\hline
			\end{tabular}
		}
		
		\vspace{0.4em}
		\begingroup
		\fontsize{9.5pt}{12pt}\selectfont
		\noindent
		$\textsuperscript{a}$\enspace``$\mathord{*}\mathord{*}\mathord{*}$'' means that $V_{\theta}f$ may not lie in $L^{2}(\mathbb{R}^{2})$ and the comparison
		is not well defined.
		\par
		\endgroup
	\end{table}

	\subsection{Structure}
	The remainder of this paper is organized as follows. In section \ref{qianyan}, we introduce the notations and collect the preliminary results needed later. Section \ref{suanzixingzhidiyibufen} is devoted to the  characterizations  of  AV-LT with a fixed opening angle. They include convolution representation, the operator norm, $L^2$-boundedness and  Fourier transform representation.  In section \ref{sobolevguanghuaruhebianhua},  
	we establish  a sharp estimate for the Sobolev smoothness of 
	$V_{\theta,\mu}f$  where $f\in L^{2}(\mathbb{R}^{2})$. In section \ref{zhijiruhebianhuan}, we establish a necessary and sufficient condition on compactly supported $f\in L^{2}(\mathbb{R}^{2})$ under which $V_{\theta,\mu}f$ is also compactly supported.
	In section \ref{daixian}, based on AV-LT samples we develop a  sampling  formula    for the reconstruction of  band-limited functions. The key ingredient  for developing such a   formula   is the design of reconstruction kernel. 
	We also establish characterizations for the sampling system   including  the optimal Riesz bounds and the 
	optimal convergence rate. 
	Finally, based on section \ref{daixian}, in  section \ref{Sobolev space jin si} we  establish  the 
	the approximation to the functions in Sobolev space
	$H^{s}(\mathbb{R}^{2})$ with $s>2$ from their AV-LT samples. Numerical simulations are conducted to verify the results in sections \ref{sec:numerical} and \ref{sec:numerical2}. In particular, the simulations  demonstrate that our reconstruction method is superior to the regularization method in terms of robustness to noise, when reconstructing band-limited functions.

	\section{Preliminaries}\label{qianyan}
	\subsection{Notations and definitions}
	
	Throughout the paper, $A\lesssim B$ means $A\leq c B$ for some constant
	$c>0$, $A\gtrsim B$ means $B\lesssim A$, and $A\asymp B$ means both
	$A\lesssim B$ and $B\lesssim A$.
	
	For $f\in L^2(\mathbb{R}^2)$, its Fourier transform and inverse Fourier transform are defined as 
	$\widehat{f}(\bm{\xi})
	:= \int_{\mathbb{R}^2} f(\bm{x})e^{-\mathrm{i}\bm{x}\cdot\bm{\xi}}d\bm{x}$
	and 
	$\widecheck{f}(\bm{x})
	:= \frac{1}{(2\pi)^2} \int_{\mathbb{R}^2} f(\bm{\xi})e^{\mathrm{i}\bm{\xi}\cdot\bm{x}}d\bm{\xi},$
	respectively. 
	Additionally, if $f$ has compact support,  then  $\widehat f$ can be extended 
	to $\mathbb C^{2}$ by the same formula $$\widehat f(\bm{\zeta}):=\int_{\mathbb R^{2}} f(\bm{x})e^{-\mathrm{i} \bm{x}\cdot \bm{\zeta}}d\bm{x},
	\quad \bm{\zeta}\in \mathbb C^{2}$$
	which defines an entire function, known as  the Fourier-Laplace transform of $f$. Given normed spaces $X$ and $Y$, we write $\mathcal{L}(X,Y)$ for the space of all
	bounded linear operators $T:X\to Y$, equipped with the operator norm
	$\|T\|_{\mathcal{L}(X,Y)}:=\sup_{\|x\|_X=1}\|Tx\|_Y$.  
	
	We recall some function spaces.
	For $\lambda\ge 0$, the Sobolev space $H^\lambda(\mathbb{R}^2)$ is defined by
	\[
	H^\lambda(\mathbb{R}^2)
	:=
	\left\{
	f:
	(1+\|\bm{\xi}\|_{2}^2)^{\lambda/2}\widehat{f}(\bm{\xi})\in L^2(\mathbb{R}^2)
	\right\},
	\]
	where $\|f\|_{H^\lambda(\mathbb R^2)}^2
	:=
	\frac{1}{(2\pi)^2}
	\int_{\mathbb R^2}
	(1+\|\bm{\xi}\|_2^2)^\lambda |\widehat f(\bm{\xi})|^2d\bm{\xi}$ (cf. \cite{Sobolevsmoothness}). 
	We denote by $H_c^\lambda(\mathbb R^2)$ the space of compactly supported elements of $H^\lambda(\mathbb R^2)$, and by $C_c^\infty(\mathbb R^2)$  the space of infinitely differentiable and  compactly supported functions on $\mathbb R^2$. 
	


	

	
	\subsection{Shannon-Whittaker sampling theorem}
	Some notation and definitions are helpful for
	introducing the Shannon-Whittaker sampling theorem. 
	Let $\Omega=[a_1,b_1]\times [a_2,b_2]\subseteq \mathbb R^2.$ Associated with  $\Omega$, we
	set $\boldsymbol{h}=(h_1,h_2)$ and  $\mathbf{c}=(c_1, c_2)$, where
	\begin{align}\label{canshu123}
		h_j:=\frac{2\pi}{b_j-a_j}, \quad  c_j:=\frac{a_j+b_j}{2}.
	\end{align}
	For $\boldsymbol{k}=(k_1,k_2)\in\mathbb Z^2$
	and $\boldsymbol{x}=(x_{1}, x_{2})$, throughout the paper  we use the notation
	\begin{align}\label{caiyangsanshu}
		\lambda_{\boldsymbol{k}}:=(h_1k_1,h_2k_2), \ \ 
		\boldsymbol{x}\oslash \boldsymbol{h}:=\left(\frac{x_1}{h_1},\frac{x_2}{h_2}\right),
	\end{align}
	where $\oslash$ denotes the Hadamard division.
	
	The Paley--Wiener space associated with $\Omega$ is defined by
	\begin{align}\label{daixianxinhaokongjian}
		PW_\Omega
		:=
		\left\{
		f\in L^2(\mathbb R^2):
		\operatorname{supp}\widehat f\subseteq \Omega
		\right\}.
	\end{align}
	The  sinc function is defined as 
	\[
	\operatorname{sinc}(t)
	:=
	\begin{cases}
		\dfrac{\sin \pi t}{\pi t}, & t\neq 0,\\[6pt]
		1, & t=0.
	\end{cases}
	\]
	Associated with sinc,  the reconstruction kernel  $S$
	for $PW_\Omega$ is defined by 
	\[
	\mathcal S(\mathbf{x}):=\operatorname{sinc}(x_1)\operatorname{sinc}(x_2), \  \forall\mathbf{x}=(x_1,x_2)\in\mathbb R^2.  
	\]
	The famous Shannon-Whittaker sampling theorem (cf. \cite[Chapter III, Theorem 1.1]{Natterer2001}) is  stated as follows.
	
	\begin{lem}\label{lem:Shannon}
		Let $f\in PW_\Omega$. Then 
		\[
		f(\mathbf{x})
		=
		\sum_{\mathbf{k}\in\mathbb Z^2}
		f(\lambda_{\boldsymbol{k}})
		e^{\mathrm{i} \mathbf{c}\cdot (\mathbf{x}-\lambda_{\boldsymbol{k}})}
		\mathcal S\left(\boldsymbol{x}\oslash \boldsymbol{h}-\mathbf{k}\right),
		\]
		where $\lambda_{\boldsymbol{k}}$ is as in \eqref{caiyangsanshu}, 
		and the  series converges  in $L^2(\mathbb R^2)$
		and pointwise on $\mathbb R^2$.
	\end{lem}
	
	\subsection{Riesz bases}
	Let $\Theta$ be a countable index set.
	A sequence $(f_k)_{k\in \Theta}$ in a Hilbert space $\mathcal H$ forms a Riesz  
	basis for $\mathcal H$ if and only if it is complete in $\mathcal H$ and satisfies the two-sided estimate
	\begin{align}\label{Rieszbudengshi}
		A\sum_k |c_k|^2
		\le
		\left\|\sum_k c_k f_k\right\|_{2}^2
		\le
		B\sum_k |c_k|^2
	\end{align}
	for some constants $A,B>0$ and all finite scalar sequences $\{c_k\}$ (cf. \cite[Theorem~1.1, $(\mathcal R_4)$]{Stoeva2020}). The constants $A$ and $B$ are called lower
	and upper Riesz  bounds, respectively. 
	The optimal Riesz lower and upper bounds are denoted by $A_{\mathrm{opt}}$ and $B_{\mathrm{opt}}$, respectively. They are important quantities since 
	the condition number (cf.\cite{Christensen,NYoung})
	of  $(f_k)_{k\in\Theta}$ is defined as $\kappa:=\sqrt{\frac{B_{\mathrm{opt}}}{A_{\mathrm{opt}}}}$.   Additionally, they are intrinsically interesting from the perspective of mathematical theory (cf.\cite{Cahill2013,Fickus2016,Massey2014}).

	
	\section{Characterization of attenuated V-line transform}\label{kehua}
	\subsection{\texorpdfstring{$L^{2}$}{}-boundedness and operator norm of AV-LT: answers to \textbf{Q1} and \textbf{Q2}}
	\label{suanzixingzhidiyibufen}
	The main purpose of this subsection is to establish the $L^{2}$-boundedness and give the  operator norm of AV-LT. 
	This requires  the  convolution and Fourier multiplier representations    of AV-LT, which   are interesting on their own.  
	We start with  Young's inequality for convolutions with finite measures.
	
	\begin{lem}(cf. \cite[Theorem 1.3]{SteinWeiss})\label{thm:young-measure}
		Let $\sigma$ be a finite complex Borel measure on $\mathbb{R}^2$.  Denote by
		$|\sigma|$ its total variation. For  
		$f\in L^2(\mathbb{R}^2)$,  define 
		$$
		(f*\sigma)(\mathbf{x})
		:= \int_{\mathbb{R}^2} f(\mathbf{x}-\mathbf{y})d\sigma(\mathbf{y}),  \ \mathbf{x}\in\mathbb{R}^2,$$
		whenever the integral exists. Then $f*\sigma$ is well defined for a.e.\ $\mathbf{x}$,
		belongs to $L^2(\mathbb{R}^2)$, and $\|f*\sigma\|_{L^2(\mathbb{R}^2)}
		\le \|f\|_{L^2(\mathbb{R}^2)}\|\sigma\|_{\mathrm{TV}}$ with $\|\sigma\|_{\mathrm{TV}}:= |\sigma|(\mathbb{R}^2)$.
	\end{lem}
	
	Next we define a Borel measure, which  will be helpful for the convoluton representation of 
	AV-LT.
	
	\begin{prop}\label{def:measure}
		Let $\mu>0$ and $\theta\in(0,\tfrac{\pi}{2})$. For every Borel set
		$E\subset\mathbb{R}^2$,  define  
		\[
		\sigma_{\mu,\theta}(E)
		:= \sum_{s\in\{+,-\}}
		\int_0^\infty \chi_E\bigl(-r \mathbf{d}_s(\theta)\bigr)
		e^{-\mu r}dr,
		\]
		where $\chi_E$ is the characteristic function of $E$,  and $\mathbf{d}_s(\theta)$
		is as in \eqref{Vline}. 
		Then $\sigma_{\mu,\theta}$ is a Borel measure on $\mathbb{R}^2$.  
	\end{prop}
	\begin{proof}
		The proof is provided in Appendix section \ref{proofofpropositionmearue}. 
	\end{proof}
	
	\begin{rem}\label{Rshangyouxiancedu}
		A direct computation shows that 
		\[
		\sigma_{\mu,\theta}(\mathbb{R}^2)
		= \sum_{s\in\{+,-\}}\int_0^\infty \chi_{\mathbb{R}^2}(-r\bm d_s(\theta))e^{-\mu r}dr
		= 2\int_0^\infty e^{-\mu r}dr
		= \frac{2}{\mu}<\infty.
		\]
		Thus $\sigma_{\mu,\theta}$ is a finite Borel measure with $\|\sigma_{\mu,\theta}\|_{\mathrm{TV}}:=|\sigma_{\mu,\theta}|(\mathbb{R}^2)=\sigma_{\mu,\theta}(\mathbb{R}^2)=\frac{2}{\mu}.$
	\end{rem}
	
	We will use the formula for change of variables in the following   form.
	
	\begin{lem}(cf.\cite[Chapter VIII, Theorem C]{GTM18})\label{thm:change-of-variables}
		Let $(X,\Sigma_X,\lambda)$ be a measure space while  $(Y,\Sigma_Y)$ a measurable
		space. Moreover, a map $T:X\to Y$ is  measurable.
		Define the measure  $\nu:=\lambda\circ T^{-1}$ by
		$\nu(B):= \lambda\bigl(T^{-1}(B)\bigr),  B\in\Sigma_Y.$
		Then for every $\nu$-integrable $f:Y\to\mathbb{C}$, we have 
		\[
		\int_Y f(y)d\nu(y)
		= \int_X f\bigl(T(x)\bigr)d\lambda(x).
		\]
	\end{lem}
	It is ready to establish the convolution representation for AV-LT as follows.
	
	\begin{lem}[Convolution representation]\label{prop:convolution}
		Let $\mu>0$, $\theta\in(0,\tfrac{\pi}{2})$ and $f\in L^2(\mathbb{R}^2)$. Then
		\[
		(V_{\theta,\mu} f)(\mathbf{x})
		= (f * \sigma_{\mu,\theta})(\mathbf{x}),
		\ \  \mathbf{x}\in\mathbb{R}^2.
		\]
	\end{lem}

	\begin{proof}
		Using Proposition \ref{def:measure},
		\begin{align*}
			(f * \sigma_{\mu,\theta})(\mathbf{x})
			&= \int_{\mathbb{R}^2}
			f(\mathbf{x}-\mathbf{y})\,d\sigma_{\mu,\theta}(\mathbf{y}) \\
			&= \sum_{s\in\{+,-\}} \int_0^\infty
			f\bigl(\mathbf{x}-(-r\mathbf{d}_s(\theta))\bigr)e^{-\mu r}\,dr \\
			&= \sum_{s\in\{+,-\}} \int_0^\infty
			f\bigl(\mathbf{x}+r\mathbf{d}_s(\theta)\bigr)e^{-\mu r}\,dr \\
			&= (V_{\theta,\mu}f)(\mathbf{x}).
		\end{align*}
		where the change of variables is ensured  by Lemma \ref{thm:change-of-variables}.
	\end{proof}
	
	By Remark \ref{Rshangyouxiancedu}, the  Borel measure $\sigma_{\mu,\theta}$ on $\mathbb{R}^2$ is finite. 
	We define its Fourier
	transform by
	\[
	\widehat{\sigma}_{\mu,\theta}(\boldsymbol{\xi})
	:= \int_{\mathbb{R}^2} e^{-\mathrm{i}\mathbf{x}\cdot\boldsymbol{\xi}}d\sigma_{\mu,\theta}(\mathbf{x}),
	\ \  \boldsymbol{\xi}\in\mathbb{R}^2.
	\]
	It has an explicit formula as follows. 
	\begin{lem}\label{lem:sigma-hat}
		For every $\boldsymbol{\xi}\in\mathbb{R}^2$, we have 
		\[
		\widehat{\sigma}_{\mu,\theta}(\boldsymbol{\xi})
		= \frac{1}{\mu - \mathrm{i}\mathbf{d}_+(\theta)\cdot\boldsymbol{\xi}}
		+ \frac{1}{\mu - \mathrm{i}\mathbf{d}_-(\theta)\cdot\boldsymbol{\xi}}.
		\]
	\end{lem}
	
	\begin{proof}
		By Proposition \ref{def:measure},
		\begin{align*}
			\widehat{\sigma}_{\mu,\theta}(\boldsymbol{\xi})
			&= \int_{\mathbb{R}^2} e^{-\mathrm{i}\mathbf{x}\cdot\boldsymbol{\xi}}d\sigma_{\mu,\theta}(\mathbf{x})
			= \sum_{s\in\{+,-\}} \int_0^\infty
			e^{-\mathrm{i}(-r\mathbf{d}_s(\theta))\cdot\boldsymbol{\xi}}e^{-\mu r}dr \\
			&= \sum_{s\in\{+,-\}}\int_0^\infty
			e^{-(\mu - \mathrm{i}\mathbf{d}_s(\theta)\cdot\boldsymbol{\xi})r}dr
			= \sum_{s\in\{+,-\}}
			\frac{1}{\mu - \mathrm{i}\mathbf{d}_s(\theta)\cdot\boldsymbol{\xi}},
		\end{align*}
		where the convergence of the  exponential integrals is derived from $\mu>0$. 
	\end{proof}
	
	Based on  Lemma \ref{prop:convolution} and 
	Lemma \ref{lem:sigma-hat}, we represent  $\widehat{V_{\theta,\mu} f}(\boldsymbol{\xi})$ through 
	a Fourier multiplier.
	
	\begin{lem}\label{Fourier-multiplier}
		Let $f\in L^2(\mathbb{R}^2)$ and $V_{\theta,\mu}f$ be as in
		\eqref{eq:Def-V-line-transform}. Then
		\begin{equation}\label{eq:Fourier identity}
			\widehat{V_{\theta,\mu} f}(\boldsymbol{\xi})
			= m_{\theta,\mu}(\boldsymbol{\xi})\widehat{f}(\boldsymbol{\xi}),
			\quad \boldsymbol{\xi}\in\mathbb{R}^2,
		\end{equation}
		where the Fourier multiplier 
		\begin{equation}\label{eq:m-theta-mu}
			m_{\theta,\mu}(\boldsymbol{\xi})
			:=\widehat{\sigma}_{\mu,\theta}(\boldsymbol{\xi})= \frac{1}{\mu - \mathrm{i}\mathbf{d}_+(\theta)\cdot\boldsymbol{\xi}}
			+ \frac{1}{\mu - \mathrm{i}\mathbf{d}_-(\theta)\cdot\boldsymbol{\xi}}.
		\end{equation}
	\end{lem}
	
	\begin{rem}\label{rem:m-properties}
		For all $\boldsymbol{\xi}\in\mathbb{R}^2$, $\bigl|\mu - \mathrm{i}\mathbf{d}_\pm(\theta)\cdot\boldsymbol{\xi}\bigr|\ge\mu>0.$
		From this, 
		one can check that 
		$m_{\theta,\mu}\in C^\infty(\mathbb{R}^2)$.
		On the other hand, 
		$\Re[m_{\theta,\mu}(\boldsymbol{\xi})]>0$ for all $\boldsymbol{\xi}\in\mathbb{R}^2$. Therefore, $1/m_{\theta,\mu}$ is also smooth on $\mathbb{R}^2$.
	\end{rem}
	
	Next we give the answers to $\textbf{Q}1$ and $\textbf{Q}2$ in section \ref{somemotivations}.

	\begin{theo}\label{L2-boundedness}
		Let $\mu>0$ and $0<\theta<\pi/2$ be fixed. Then  $V_{\theta,\mu}\in \mathcal{L}\bigl(L^2(\mathbb R^2), L^2(\mathbb R^2)\bigr)$   and $\|V_{\theta,\mu}\|_{L^2(\mathbb R^2)\to L^2(\mathbb R^2)}=2/\mu$.
		
	\end{theo}
	
	\begin{proof}
		The linearity of $V_{\theta,\mu}$ is immediate from Definition \ref{eq:Def-V-line-transform}. By
		Lemma \ref{prop:convolution}, we have 
		$V_{\theta,\mu}f=f*\sigma_{\mu,\theta}$. By Remark \ref{Rshangyouxiancedu}, 
		$\|\sigma_{\mu,\theta}\|_{\scriptscriptstyle\mathrm{TV}}=2/\mu$. Then  Young's inequality for convolution with a
		finite measure gives $\|V_{\theta,\mu}f\|_{L^2(\mathbb R^2)}\leq (2/\mu)\|f\|_{L^2(\mathbb R^2)}$ for all
		$f\in L^2(\mathbb R^2)$.
		Thus $\|V_{\theta,\mu}\|_{L^2(\mathbb R^2)\to L^2(\mathbb R^2)}\leq 2/\mu<\infty$. 
		
		It remains to show that the  upper bound $2/\mu$ is sharp. For every nonzero
		$f\in L^2(\mathbb R^2)$, Theorem \ref{Fourier-multiplier} and Plancherel's theorem give
		\[
		\frac{\|V_{\theta,\mu}f\|_{L^2(\mathbb{R}^2)}}{\|f\|_{L^2(\mathbb{R}^2)}}
		=
		\frac{\|m_{\theta,\mu}\widehat f\|_{L^2(\mathbb{R}^2)}}{\|\widehat f\|_{L^2(\mathbb{R}^2)}}.
		\]
		By Remark \ref{rem:m-properties}, the multiplier
		$m_{\theta,\mu}$ is continuous and satisfies $m_{\theta,\mu}(\boldsymbol{0})=2/\mu$.
		Hence, for every $\epsilon>0$, there exists $r>0$ such that
		$|m_{\theta,\mu}(\boldsymbol{\xi})|\geq 2/\mu-\epsilon$ for all
		$\boldsymbol{\xi}\in B(\boldsymbol{0},r)$, where $B(\boldsymbol{0},r)$ denotes the open ball with radius $r$ in $\mathbb R^2$.
		
		Let $0\neq g_r\in C_c^\infty(B(\boldsymbol{0},r))$ and
		$f_r=\widecheck{g_r}$. Then $f_r\in L^2(\mathbb R^2)$ and
		$\widehat f_r=g_r$. Therefore, 
		\[
		\frac{\|V_{\theta,\mu}f_r\|_{L^2(\mathbb{R}^2)}}{\|f_r\|_{L^2(\mathbb{R}^2)}}
		=
		\frac{\|m_{\theta,\mu}g_r\|_{L^2(\mathbb{R}^2)}}{\|g_r\|_{L^2(\mathbb{R}^2)}}
		\geq
		2/\mu-\epsilon.
		\]
		Then 
		$\|V_{\theta,\mu}\|_{L^2(\mathbb{R}^2)\to L^2(\mathbb{R}^2)}\geq 2/\mu-\epsilon$. Letting
		$\epsilon\rightarrow0$, we obtain
		$\|V_{\theta,\mu}\|_{L^2(\mathbb{R}^2)\to L^2(\mathbb{R}^2)}\geq 2/\mu$. This completes the proof.
	\end{proof}

	
	The following states that $V_{\theta}$ is not well defined 
	on $L^{2}(\mathbb{R}^{2})$, which is essentially different from $V_{\mu,\theta}$.

	\begin{prop}\label{Vlinebushoulian}
		There exists $f\in L^2(\mathbb{R}^2)$ such that
		$V_{\theta}f\notin L^2(\mathbb R^2)$. That is,    $V_{\theta}$ does not map $L^2(\mathbb R^2)$ into itself.
	\end{prop}

	\begin{proof}
		We choose nonzero, nonnegative functions $\varphi,\psi\in C_c^\infty((-1,1))$ and define 
		\begin{align}\label{feijindef}
			f(\mathbf{x}):=\varphi(\mathbf{x}\cdot \mathbf{d}_+(\theta))\psi(\mathbf{x}\cdot \mathbf{d}_+^\perp(\theta)),
		\end{align}
		where $\mathbf{d}_+^\perp(\theta)$ denotes a unit vector orthogonal to $\mathbf{d}_+(\theta)$.
		Then $f\in C_c^\infty(\mathbb R^2)\subseteq
		L^2(\mathbb R^2)$. Write
		$\mathbf{x}=a\mathbf{d}_+(\theta)+b\mathbf{d}_+^\perp(\theta)$. Since
		$f\ge0$,  then 
		\[
		V_{\theta}f(\mathbf{x})
		\ge
		\int_0^\infty f(\mathbf{x}+r\mathbf{d}_+(\theta))dr
		=
		\psi(b)\int_0^\infty \varphi(a+r)dr .
		\]
		If $a\le -1$, then $\operatorname{supp}\varphi\subseteq(-1,1)$ implies
		$\int_0^\infty \varphi(a+r)dr=\int_{\mathbb R}\varphi(t)dt$. Therefore,
		for all $a\le -1$,
		\begin{align}\label{feijinzhicheng}
			V_{\theta}f(a\mathbf{d}_+(\theta)+b\mathbf{d}_+^\perp(\theta))
			\ge
			\psi(b)\int_{\mathbb R}\varphi(t)dt .
		\end{align}
		Since $\varphi$ and $\psi$ are nonzero and nonnegative, the right-hand side is
		not identically zero as a function of $b$. Hence
		\[
		\|V_{\theta}f\|_{L^2(\mathbb R^2)}^2
		\ge
		\int_{-\infty}^{-1}\int_{\mathbb R}
		\left|\psi(b)\int_{\mathbb R}\varphi(t)dt\right|^2dbda
		=\infty .
		\]
		This completes the proof. 
	\end{proof}
	
	\subsection{Sharp estimation of the Sobolev smoothness of \texorpdfstring{$V_{\theta, \mu} f$}{} and the answer to \textbf{Q3}}
	\label{sobolevguanghuaruhebianhua}
	We first establish the sharp estimation for  $\nu_{2}(V_{\theta, \mu}f)$. Herein the sharpness means the estimation
	holds for all  $f\in H^{\lambda}(\mathbb{R}^2)$
	where  $\lambda>0 $ is arbitrarily fixed. 
	After this, we give an  answer to \textbf{Q3} in section \ref{somemotivations}.

	\begin{lem}\label{lem:modulation-estimate}
		For $\lambda\ge0$ and $0 \not\equiv\psi\in  C_c^\infty(\mathbb R^2)$, we define $\psi_N(\mathbf x):=e^{\mathrm iN\mathbf e\cdot\mathbf x}\psi(\mathbf x)$ for each $N\ge0$, where $\mathbf e\in \mathbb{S}^1$. Then $\|\psi_N\|_{H^\lambda(\mathbb R^2)}\asymp(1+N^2)^{\lambda/2}$ for all $N\ge0$, where the implicit constants depend only on $\lambda$ and $\psi$.
	\end{lem}
	
	\begin{proof}
		See Appendix \ref{Sobolev-nonsmoothing} for the  proof.
	\end{proof}

	\begin{theo}\label{Vline-sobolev}
		Let $\mu>0$ and $0<\theta<\pi/2$ be fixed. 
		Then the following statements hold.
		\begin{enumerate}[label=(\roman*)]
			\item For every $\lambda\ge0$, $V_{\theta,\mu}\in \mathcal{L}\bigl(H^\lambda(\mathbb{R}^2),H^{\lambda}(\mathbb{R}^2)\bigr)$.
			\item For every $\lambda\ge0$ and every $\epsilon>0$,    
			there exists $f\in H^\lambda_c(\mathbb R^2)\subseteq
			H^\lambda(\mathbb R^2)$
			such that $V_{\theta,\mu}f\notin H^{\lambda+\epsilon}(\mathbb R^2)$.  
		\end{enumerate}   	
	\end{theo}
	
	\begin{proof}
		\underline{\textbf{Proof of (i)}}. 
		By Lemma \ref{Fourier-multiplier} and Remark \ref{rem:m-properties}, for every $f\in H^\lambda(\mathbb{R}^2)$ we have
		\begin{align}
			\|V_{\theta,\mu}f\|_{H^\lambda(\mathbb{R}^2)}^2\notag
			&=\frac1{(2\pi)^2}\int_{\mathbb R^2}(1+\left\|\boldsymbol{\xi}\right\|_2^2)^\lambda
			|m_{\theta,\mu}(\boldsymbol{\xi})|^2|\widehat f(\boldsymbol{\xi})|^2d\boldsymbol{\xi}\notag\\
			&\le \|m_{\theta,\mu}\|_{L^\infty(\mathbb R^2)}^2
			\frac1{(2\pi)^2}\int_{\mathbb R^2}(1+\left\|\boldsymbol{\xi}\right\|_2^2)^\lambda
			|\widehat f(\boldsymbol{\xi})|^2d\boldsymbol{\xi}\notag\\
			&=\frac{4}{\mu^2}\|f\|_{H^\lambda(\mathbb R^2)}^{2},\notag
		\end{align}
		from which  it follows that
		$V_{\theta,\mu}:H^\lambda(\mathbb R^2)\to H^\lambda(\mathbb R^2)$ is bounded.
		
		\underline{\textbf{Proof of (ii)}}. For fixed $\lambda\ge0$ and $\epsilon>0$, we define $$f_N(\mathbf x):=(1+N^2)^{-\lambda/2}e^{\mathrm iN\mathbf e\cdot\mathbf x}\psi(\mathbf x)$$ for every integer $N\ge 2$, where $\psi\in C_c^\infty(\mathbb R^2)$ is nonzero and nonnegative, and $\mathbf e:=(\cos\theta,-\sin\theta)$ satisfies $\mathbf e\cdot\mathbf d_+=0$. 
		It follows from Lemma \ref{lem:modulation-estimate} that $\sup_{N\geq 2}\|f_N\|_{H^\lambda(\mathbb R^2)}<\infty$.

		The direct calculation leads to  $\widehat f_N(\boldsymbol\xi)=(1+N^2)^{-\lambda/2}
		\widehat\psi(\boldsymbol\xi-N\mathbf e)$.
		For $\|\boldsymbol\eta\|_2\le1$, we have
		$1+\|N\mathbf e+\boldsymbol\eta\|_2^2\gtrsim1+N^2$,
		and \eqref{eq:m-theta-mu} gives
		\begin{equation}\label{eq:modulated-multiplier-lower}
			|m_{\theta,\mu}(N\mathbf e+\boldsymbol\eta)|
			\ge \Re[m_{\theta,\mu}(N\mathbf e+\boldsymbol\eta)]
			\ge \frac{\mu}{\mu^2+(\boldsymbol\eta\cdot\mathbf d_+)^2}
			\ge \frac{\mu}{\mu^2+1}.
		\end{equation}
		Then it follows from  \eqref{eq:Fourier identity}  that 
		\begin{align}
			\|V_{\theta,\mu}f_N\|_{H^{\lambda+\epsilon}(\mathbb{R}^{2})}^2
			&=\frac1{(2\pi)^2}
			\int_{\mathbb R^2}(1+\|\boldsymbol\xi\|_2^2)^{\lambda+\epsilon}
			|m_{\theta,\mu}(\boldsymbol\xi)|^2
			|\widehat f_N(\boldsymbol\xi)|^2d\boldsymbol\xi\notag\\
			&=\frac{(1+N^2)^{-\lambda}}{(2\pi)^2}
			\int_{\mathbb R^2}(1+\|N\mathbf e+\boldsymbol\eta\|_2^2)^{\lambda+\epsilon}
			|m_{\theta,\mu}(N\mathbf e+\boldsymbol\eta)|^2
			|\widehat\psi(\boldsymbol\eta)|^2d\boldsymbol\eta\notag\\
			&\ge\frac{\mu^2(1+N^2)^{-\lambda}}{(2\pi)^2(\mu^2+1)^2}
			\int_{\|\boldsymbol\eta\|_2\le1}
			(1+\|N\mathbf e+\boldsymbol\eta\|_2^2)^{\lambda+\epsilon}
			|\widehat\psi(\boldsymbol\eta)|^2d\boldsymbol\eta\notag\\
			&\gtrsim (1+N^2)^\epsilon
			\int_{\|\boldsymbol\eta\|_2\le1}
			|\widehat\psi(\boldsymbol\eta)|^2d\boldsymbol\eta.\label{you xia jie}
		\end{align}
		The integral is strictly positive
		because $\widehat\psi(\textbf{0})=\int_{\mathbb R^2}\psi(\mathbf x)d\mathbf x>0$ and $\widehat\psi$
		is continuous at $\textbf{0}$.
		
		Assume that $V_{\theta,\mu}f\in H^{\lambda+\epsilon}(\mathbb R^2)$
		for every $f\in H^\lambda_c(\mathbb R^2)$. With $K:=\operatorname{supp}\psi$, the space $H^\lambda_K:=\{f\in H^\lambda(\mathbb R^2):\operatorname{supp}f\subseteq K\}$ is a Banach space.
		We first   prove  that  $V_{\theta,\mu}$ has a closed graph. Actually,  suppose that $g_j\to g$ in
		$H^\lambda_K$ and $V_{\theta,\mu}g_j\to h$ in $H^{\lambda+\epsilon}(\mathbb R^2)$.
		By (i), $V_{\theta,\mu}g_j\to V_{\theta,\mu}g$ in $H^\lambda(\mathbb R^2)$.
		On the other hand, the continuous embedding
		$H^{\lambda+\epsilon}(\mathbb R^2)\hookrightarrow H^\lambda(\mathbb R^2)$
		implies that $V_{\theta,\mu}g_j\to h$ in $H^\lambda(\mathbb R^2)$.
		From the uniqueness of limit, we get  $h=V_{\theta,\mu}g$. Then  the graph is closed.
		By the closed graph theorem, there exists a constant $C>0$ such that
		\begin{align}\label{quanjuguji}
			\|V_{\theta,\mu}f\|_{H^{\lambda+\epsilon}(\mathbb{R}^{2})}\le C\|f\|_{H^\lambda(\mathbb{R}^{2})}   
		\end{align}
		for every $f\in H^\lambda_K$.
		Since $f_N\in H^\lambda_K$ and $\sup_N\|f_N\|_{H^\lambda(\mathbb{R}^{2})}<\infty$, it follows from 
		\eqref{quanjuguji} that 
		$\sup_N\|V_{\theta,\mu}f_N\|_{H^{\lambda+\epsilon}(\mathbb{R}^{2})}\leq C \sup_N\|f_N\|_{H^\lambda(\mathbb{R}^{2})}<\infty$. This 
		contradicts \eqref{you xia jie}, and  the proof of  (ii)
		is completed.
	\end{proof}
	
	\begin{rem}
		(1) From Theorem \ref{Vline-sobolev} (ii), the  estimation $\nu_{2}(V_{\theta,\mu}f)\geq \nu_{2}(f)$ is sharp 
		for all functions $f\in H^{\lambda}(\mathbb{R}^{2})$ where $\lambda>0$ is arbitrarily fixed. 
		(2) For any $f\in L^2(\mathbb{R}^2)$, if $\nu_{2}(f)>1$
		then $\nu_{2}(V_{\theta\mu}f)>1$. This sharp  condition gives an answer  to \textbf{Q3}. 
	\end{rem}

	Theorem \ref{Vline-sobolev} shows that, for every $\epsilon>0$, the inequality
	$\nu_2(V_{\theta,\mu}f)\geq\nu_2(f)+\epsilon$ cannot hold for all $f$. However, the following example shows
	there exists $f\in L^{2}(\mathbb{R}^{2})$ such that 
	$\nu_2(V_{\theta,\mu}f)>\nu_2(f)$.
	
	\begin{exam}\label{ex:strict-sobolev}
		For fixed $\mu>0$ and $0<\theta<\pi/2$, define
		$f(x_1,x_2)=[\chi_{(0,1]}*\chi_{(0,1]}](x_1)
		[\chi_{(0,1]}*\chi_{(0,1]}](x_2)$, where $\chi_{(0,1]}$ is the characteristic
		function of the interval $(0,1]$. From \cite[Example 1.3]{HanShen},
		we have $\nu_2(f)=3/2$. In Appendix \ref{app:sobolev-details}, we prove that $\nu_2(V_{\theta,\mu}f)=5/2$.
	\end{exam}
	
	The following example  shows that the Radon transform does not necessarily improve Sobolev smoothness. This is similar as AV-LT.

	\begin{exam}\label{example3.2}
		Consider the box-spline  $B(x_{1}, x_{2})=\chi_{(0, 1]}(x_{1})\chi_{(0, 1]}(x_{2})$. Its Fourier transform  is given by
		\begin{equation}\label{eq:box-spline-fourier}
			\widehat{B}(\boldsymbol{\xi})
			=
			\frac{1-e^{-\mathrm{i}\boldsymbol{k}_1\cdot\boldsymbol{\xi}}}
			{\mathrm{i}\boldsymbol{k}_1\cdot\boldsymbol{\xi}}
			\frac{1-e^{-\mathrm{i}\boldsymbol{k}_2\cdot\boldsymbol{\xi}}}
			{\mathrm{i}\boldsymbol{k}_2\cdot\boldsymbol{\xi}},
			\quad \boldsymbol{\xi}\in\mathbb{R}^2,
		\end{equation}
		where $\boldsymbol{k}_1=(1,0)$ and $\boldsymbol{k}_2=(0,1)$.
		Taking $\vartheta=(0,1)$ in \eqref{Radonbianhuandedingyi}, from the  Fourier slice theorem (cf. \cite{Han}) we obtain $\widehat{R_{\vartheta}B}(\xi)=\frac{1-e^{-\mathrm{i}\xi}}{\mathrm{i}\xi}, \  \xi\in \mathbb{R}$. By \cite[Page 375]{HanShen}, we have $\nu_{2}(B)=\nu_{2}(R_{\vartheta}B)=1/2$. 
	\end{exam}
	
	\subsection{Answer to \textbf{Q}4: characterization of compact support for attenuated V-line transform}
	\label{zhijiruhebianhuan}
	The purpose of this subsection is to establish a sufficient and  necessary   condition on  the  compactly supported function $f\in L^2(\mathbb R^2)$ such that 
	$V_{\theta,\mu}f$ is also compactly supported.
	
	\begin{theo}\label{thm:compact-support-V}
		Define $\boldsymbol d_+^\perp(\theta):=(\cos\theta,-\sin\theta)$ and
		$\boldsymbol d_-^\perp(\theta):=(\cos\theta,\sin\theta)$. Consequently, 
		$\boldsymbol d_\pm(\theta)\cdot\boldsymbol d_\pm^\perp(\theta)=0$.
		Let $\mu>0$ and assume that
		$f\in L^2(\mathbb R^2)$ is compactly supported. Then the following two  items   are equivalent:
		\begin{enumerate}
			\item[(1)] $V_{\theta,\mu}f$ is compactly supported.
			\item[(2)] The Fourier-Laplace transform of $f$ vanishes on each of the two   complex lines
			$\{\omega\boldsymbol d_\pm^\perp(\theta)-\mathrm i\mu\boldsymbol d_\pm(\theta):\omega\in\mathbb R\}\subseteq\mathbb C^2$.
		\end{enumerate}
	\end{theo}
	
	\begin{proof}
		Let $K$ be a convex compact set such that $\operatorname{supp}f\subseteq K\subseteq\mathbb R^2$. Since $f\in L^2(\mathbb R^2)$ has compact support, it follows from  Fubini's theorem that the two  integrals
		\begin{equation}\label{eq:weighted-line-integral}
			I_\pm f(r):=\int_{\mathbb R}
			f\bigl(r\boldsymbol d_\pm^\perp(\theta)+u\boldsymbol d_\pm(\theta)\bigr)
			e^{-\mu u}du
		\end{equation}
		are well defined for a.e. $r\in \mathbb R$. Moreover, both $I_+f(r)$ and $I_-f(r)$ belong to $L^1(\mathbb R)$. On the other hand, through the change of variables one gets 
		\begin{equation}\label{eq:weighted-slice}
			\widehat{I_\pm f}(\omega)
			=\widehat f\bigl(\omega\boldsymbol d_\pm^\perp(\theta)
			-\mathrm i\mu\boldsymbol d_\pm(\theta)\bigr),
			\qquad \omega\in\mathbb R.
		\end{equation}
		
		\underline{\bm{$(1)\Longrightarrow(2)$}}.
		By \eqref{eq:weighted-slice}, 
		it suffices to prove that 
		\begin{align}
			I_+f(r)=I_-f(r)=0  \ \hbox{a.e.} \  r\in \mathbb R.
		\end{align}
		We prove the assertion for $I_+f$, and  the proof for $I_-f$ is analogous.
		
		Since $K$ and $\operatorname{supp}(V_{\theta,\mu}f)$ are compact, there exists $M>0$ such that every $\boldsymbol y\in K\cup\operatorname{supp}(V_{\theta,\mu}f)$ satisfies  $\|\boldsymbol y\|_2\leq M$. 
		For $r,a\in\mathbb R$, set 
		\begin{align}\label{xheArdingyi}
			\boldsymbol x:=r\boldsymbol d_+^\perp(\theta)+a\boldsymbol d_+(\theta),  A_r:=\min\left\{
			-M,-\frac{M+r\cos(2\theta)}{\sin(2\theta)}
			\right\}.  
		\end{align}
		The proof proceeds in  the following   Steps 1-3 and applies to almost every
		$(r,a)\in\Lambda$, where
		\[
		\Lambda:=\{(r,a)\in\mathbb R^2:a<A_r\},
		\]
		and 
		\newlength{\steplabelwidth}
		\settowidth{\steplabelwidth}{\textbf{Step 3:}\quad}
		
		\noindent\hspace*{1em}%
		\makebox[\steplabelwidth][l]{\textbf{Step 1:}}%
		We prove $(V_{\theta,\mu}f)(\boldsymbol x)=0$ and
		$(T_+f)(\boldsymbol x)=e^{\mu a}I_+f(r)$;
		
		\noindent\hspace*{1em}%
		\makebox[\steplabelwidth][l]{\textbf{Step 2:}}%
		We prove $(T_-f)(\boldsymbol x)=0$;
		
		\noindent\hspace*{1em}%
		\makebox[\steplabelwidth][l]{\textbf{Step 3:}}%
		By \eqref{eq:attenuated-ray-transforms-right}
		and \eqref{eq:attenuated-ray-transforms-left}, we have
		\begin{align}\label{gongshi312}
			(V_{\theta,\mu}f)(\boldsymbol x)
			=(T_+f)(\boldsymbol x)+(T_-f)(\boldsymbol x).
		\end{align}
		
		\noindent\hspace*{1em}\hspace*{\steplabelwidth}%
		Combining Steps 1--2 and \eqref{gongshi312}, we obtain $I_+f(r)=0$.


		\underline{\textbf{Proof for Step 1:}}
		For $a<A_r\leq -M$, we have 
		\begin{align}\label{xiangliangchangduguji}
			\|\bm{x}\|_{2}=
			\|r\boldsymbol d_+^\perp(\theta)
			+a\boldsymbol d_+(\theta)\|_{2}=\sqrt{r^2+a^2}>M.
		\end{align}
		Then $(V_{\theta,\mu}f)(\boldsymbol x)=0$. Next we prove $(T_+f)(\boldsymbol x)=e^{\mu a}I_+f(r)$. 
		By the change of variables $u=a+t$ in \eqref{eq:attenuated-ray-transforms-right}, we obtain
		\begin{equation}\label{eq:ray}
			(T_+f)(\boldsymbol x)
			=e^{\mu a}\int_a^\infty
			f\bigl(r\boldsymbol d_+^\perp(\theta)+u\boldsymbol d_+(\theta)\bigr)
			e^{-\mu u}du.
		\end{equation}
		We split
		\begin{equation}\label{I_+ split}
			I_+f(r)
			=
			\underbrace{
				\int_{-\infty}^{a}
				f\bigl(r\boldsymbol d_+^\perp(\theta)
				+u\boldsymbol d_+(\theta)\bigr)e^{-\mu u}du
			}_{I_{+,1}}
			+
			\underbrace{
				\int_a^\infty
				f\bigl(r\boldsymbol d_+^\perp(\theta)
				+u\boldsymbol d_+(\theta)\bigr)e^{-\mu u}du
			}_{I_{+,2}}.
		\end{equation}
		Through the similar  analysis as  in 
		\eqref{xiangliangchangduguji}, we have  $r\boldsymbol d_+^\perp(\theta)
		+u\boldsymbol d_+(\theta)\notin K$ for $u\leq a$.
		Since $\operatorname{supp}f\subseteq K$, it follows that
		$I_{+,1}=0$. Then it follows from 
		\eqref{eq:ray} and \eqref{I_+ split} that 
		$(T_+f)(\boldsymbol x)=e^{\mu a}I_+f(r)$.
		
		\underline{\textbf{Proof for Step 2:}} 
		It is sufficient  to prove the ray  $\{\boldsymbol x+t\boldsymbol d_-(\theta):t\geq0\}$ is disjoint from $K$. 
		Actually, for  every $t\geq0$, it follows from   
		\eqref{xheArdingyi}, $a<A_r$ and  $\sin(2\theta)>0$ that 
		$$
		\bigl(\boldsymbol x+t\boldsymbol d_-(\theta)\bigr)
		\cdot\boldsymbol d_-^\perp(\theta)
		=r\cos(2\theta)+a\sin(2\theta)
		<-M.
		$$
		However,  $\boldsymbol y\cdot\boldsymbol d_-^\perp(\theta)\geq-M$ for every $\boldsymbol y\in K$. The proof is completed. 
		
		\underline{\textbf{Proof for Step 3:}} Combining \textbf{Steps 1} and 2 with \eqref{gongshi312}, we obtain,  for almost every $(r,a)\in\Lambda$, 
		$$
		0=(V_{\theta,\mu}f)(\boldsymbol x)=e^{\mu a}I_+f(r).
		$$
		Since $e^{\mu a}>0$, it follows from Fubini's theorem that $I_+f(r)=0$ for a.e. $r\in\mathbb R$. This proves (2).
		
		\underline{\bm{$(2)\Longrightarrow(1)$}}.
		By (2) and  \eqref{eq:weighted-slice}, we have $I_+f(r)=I_-f(r)=0$ for a.e. $r\in\mathbb R$. We just need to  prove $(T_+f)(\boldsymbol x)=0$ for a.e. $\boldsymbol x\notin K$, and the proof for $T_-f$ is analogous. For $\boldsymbol x\notin K$, the convexity of $K$ implies that at least one of the two opposite rays
		$\{\boldsymbol x+t\boldsymbol d_+(\theta):t\geq0\}$ and
		$\{\boldsymbol x-t\boldsymbol d_+(\theta):t\geq0\}$
		is disjoint from $K$. If both rays are disjoint from $K$, it is clear that $(T_+f)(\boldsymbol x)=0$. Therefore,  we just need to consider the case illustrated  in 
		Fig.\ref{fig:convex-rays} where only one 
		ray intersets with $K$. 
		\begin{figure}[H]
			\centering
			\includegraphics[
			width=0.6\textwidth,
			height=3cm
			]{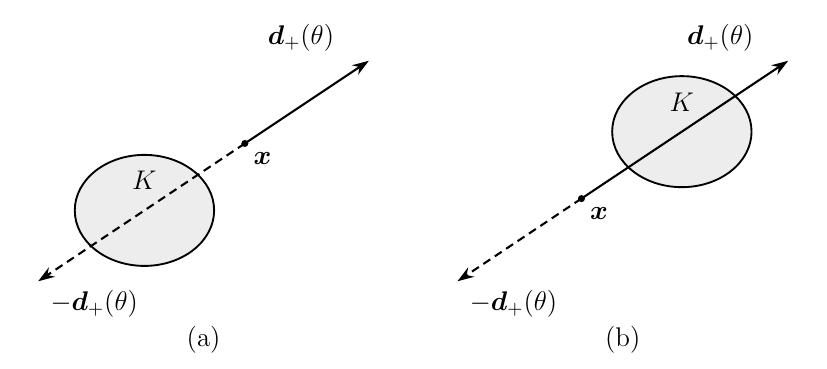}
			\caption{Rays disjoint from $K$: (a) along $\boldsymbol d_+(\theta)$; (b) along $-\boldsymbol d_+(\theta)$.}
			\label{fig:convex-rays}
		\end{figure}
		Recall that  $\boldsymbol x=r\boldsymbol d_+^\perp(\theta)+a\boldsymbol d_+(\theta)$. 
		If the ray $\{\boldsymbol x+t\boldsymbol d_+(\theta):t\ge0\}$ is disjoint from $K$, as illustrated in Fig.\ref{fig:convex-rays}(a), then \eqref{eq:ray} and $\operatorname{supp}f\subseteq K$ give  $(T_+f)(\boldsymbol x)=0$.
		
		If the ray $\{\boldsymbol x-t\boldsymbol d_+(\theta):t\ge0\}$ is disjoint from $K$, as illustrated in Fig.\ref{fig:convex-rays}(b), then for every $u\le a$, we have
		$$
		r\boldsymbol d_+^\perp(\theta)+u\boldsymbol d_+(\theta)
		=\boldsymbol x-(a-u)\boldsymbol d_+(\theta)\notin K.
		$$
		Since $\operatorname{supp}f\subseteq K$, we obtain
		$f(r\boldsymbol d_+^\perp(\theta)+u\boldsymbol d_+(\theta))=0$
		for a.e. $u\leq a$. Therefore, \eqref{eq:ray} gives 
		$$
		\begin{aligned}
			(T_+f)(\boldsymbol x)
			&=e^{\mu a}\int_a^\infty
			f\bigl(r\boldsymbol d_+^\perp(\theta)+u\boldsymbol d_+(\theta)\bigr)
			e^{-\mu u}du\\
			&=e^{\mu a}\int_{\mathbb R}
			f\bigl(r\boldsymbol d_+^\perp(\theta)+u\boldsymbol d_+(\theta)\bigr)
			e^{-\mu u}du\\
			&=e^{\mu a}I_+f(r)=0.
		\end{aligned}
		$$
		Thus $(T_+f)(\boldsymbol x)=0$ for a.e. $\boldsymbol x\notin K$. So $\operatorname{supp}(V_{\theta,\mu}f)\subseteq K$, which proves (1).
	\end{proof}

	The following constructs a class of compactly supported functions in $L^{2}(\mathbb{R}^{2})$ such that their 
	attenuated V-line transforms are also compactly supported. 
	
	\begin{prop}\label{compact-example}
		Let  $0\not\equiv g\in L^2(\mathbb R^2)$ be an arbitrary and compactly supported function. 
		Associated with fixed $\theta$ and $\mu>0$, 
		define a compactly supported function  $f:\mathbb R^2\to\mathbb C$ by
		\[
		\begin{aligned}
			f(\mathbf x)
			:=
			&-\frac14\Big(
			e^{-2\mu}g(\mathbf x+\mathbf d_{+}(\theta)+\mathbf d_{-}(\theta))
			- g(\mathbf x+\mathbf d_{+}(\theta)-\mathbf d_{-}(\theta))\\
			&\qquad\ - g(\mathbf x-\mathbf d_{+}(\theta)+\mathbf d_{-}(\theta))
			+ e^{2\mu}g(\mathbf x-\mathbf d_{+}(\theta)-\mathbf d_{-}(\theta))
			\Big),\ \mathbf x \in\mathbb{R}^{2}.
		\end{aligned}
		\]
		Then  $V_{\theta,\mu}f\in L^{2}(\mathbb{R}^{2})$ is compactly supported. 
	\end{prop}
	\begin{proof}
		The proof is provided  in Appendix \ref{app:compact-example-details}.
	\end{proof}
	
	We provide a  counterexample as follows. 
	
	\begin{theo}
		\label{compact-remark}
		For a nonzero compactly supported
		function $0\leq f\in L^2(\mathbb{R}^{2}),$
		$V_{\theta,\mu}f$ is not compactly supported. 
	\end{theo}
	\begin{proof}
		If  $V_{\theta,\mu}f$ is compactly supported,
		then it follows from Theorem \ref{thm:compact-support-V} (2) that $\widehat f\!\bigl(-\mathrm{i}\mu\mathbf d_s(\theta)\bigr)=0$ for each $s\in\{+,-\}$.
		On the other hand, by the definition of the Fourier--Laplace
		transform, we have 
		\[
		\widehat f\!\bigl(-\mathrm{i}\mu\mathbf d_s(\theta)\bigr)
		=\int_{\mathbb R^2} f(\mathbf x)e^{-\mu\mathbf x\cdot \mathbf d_s(\theta)}d\mathbf x.
		\]
		Clearly, $\widehat f(-\mathrm{i}\mu\mathbf d_s(\theta))>0$ which leads to a  contradiction to the above. The proof is completed.    
	\end{proof}

	\begin{rem}\label{bumanzujinzhcheng}
		It follows from \eqref{feijindef} and \eqref{feijinzhicheng} that 
		V-LT does not necessarily preserve compact support. On the other hand, setting $\mu=0$ in the proof of
		Theorem \ref{thm:compact-support-V}   shows that, for any compactly supported $f\in L^2(\mathbb R^2)$, $V_\theta f$ is compactly supported if and only if $\widehat f\bigl(\omega\boldsymbol d_{\pm}^{\perp}(\theta)\bigr)=0$ for all $\omega\in\mathbb R$.
	\end{rem}

	\section{ Reconstruction formula   for  band-limited functions from their  AV-LT samples and its characterization}\label{daixian}
	
	In this section, we establish a reconstruction formula for band-limited functions from their AV-LT samples. The reconstruction is implemented by directly evaluating the formula, without solving any equation.
	
	\subsection{Answer to \textbf{Q5} for bandlimited case: reconstruction formula  from AV-LT sample}
	We begin with the reconstruction kernel function constructed  as follows.
	
	\begin{defi}\label{def:PhiN}
		Let the   Paley--Wiener space  $PW_\Omega$  be  as  in \eqref{daixianxinhaokongjian}. Define
		\begin{align}\label{TTCX}
			\Phi(\boldsymbol{x})
			:=
			\frac{h_1h_2}{(2\pi)^2}
			\int_\Omega
			\frac{1}{m_{\theta,\mu}(\boldsymbol{\xi})}
			e^{\mathrm{i}\boldsymbol{x}\cdot \boldsymbol{\xi}}d\boldsymbol{\xi},
			\ \boldsymbol{x}\in\mathbb R^2,
		\end{align}
		where $h_{j}$ and  $m_{\theta,\mu}(\boldsymbol{\xi})$ are  as in \eqref{canshu123} and  \eqref{eq:m-theta-mu}, respectively. Equivalently,
		\begin{align}\label{Fourierexpression}
			\widehat{\Phi}(\boldsymbol{\xi})
			=h_1h_2
			\frac{\chi_{\Omega}(\boldsymbol{\xi})}
			{m_{\theta,\mu}(\boldsymbol{\xi})},
			\ \boldsymbol{\xi}\in\mathbb R^2.
		\end{align}
	\end{defi}
	
	Based on $\Phi$, the sampling-based reconstruction  formula is established  as follows. 
	
	\begin{theo}\label{thm:sampling-reconstruction}
		Suppose that  $f\in PW_\Omega$.  As
		in \eqref{caiyangsanshu}, let 
		$\lambda_{\boldsymbol{k}}=
		(h_1k_1,h_2k_2)$ with $\boldsymbol{k}=(k_{1}, k_{2})$. Then for every $\boldsymbol{x}\in\mathbb R^2$, it holds that
		\begin{align}\label{caiyangdinglidaixian}
			f(\boldsymbol{x})
			=
			\sum_{\mathbf{k}\in\mathbb Z^2}
			[V_{\theta,\mu}f]
			\left(\lambda_{\boldsymbol{k}}\right)
			\Phi\left(\boldsymbol{x}-\lambda_{\boldsymbol{k}}\right).
		\end{align}
	\end{theo}
	
	\begin{proof}
		By  \eqref{eq:Fourier identity} we have
		$
		\widehat{V_{\theta,\mu}f}(\boldsymbol{\xi})
		= m_{\theta,\mu}(\boldsymbol{\xi})\widehat{f}(\boldsymbol{\xi}).
		$
		From this and $\hbox{supp}(\widehat{f})\subseteq\Omega$, we have $V_{\theta,\mu}f\in PW_\Omega$. 
		Additionally, since $m_{\theta,\mu}(\boldsymbol{\xi})\neq 0$ for all
		$\boldsymbol{\xi}\in \mathbb{R}^{2}$, it follows that
		$\widehat{f}(\boldsymbol{\xi})=\frac{1}{m_{\theta,\mu}(\boldsymbol{\xi})}\widehat{V_{\theta,\mu}f}(\boldsymbol{\xi})$.
		Then the   inverse Fourier transform  gives
		\begin{equation}\label{eq:inv}
			f(\mathbf{x})
			= \frac{1}{(2\pi)^2}\int_\Omega
			\frac{1}{m_{\theta,\mu}(\boldsymbol{\xi})}
			\widehat{V_{\theta,\mu}f}(\boldsymbol{\xi})e^{\mathrm{i}\mathbf{x}\cdot\boldsymbol{\xi}}d\boldsymbol{\xi},
			\quad a. e. \  \mathbf{x}\in\mathbb{R}^2.
		\end{equation}
		Since $\Omega$ is bounded and $m_{\theta,\mu}$ is continuous with $m_{\theta,\mu}\neq 0$ on $\Omega$, we have  $1/m_{\theta,\mu}\in L^\infty(\Omega)$.
		Moreover, it follows from $\widehat{V_{\theta,\mu}f}\in L^2(\Omega)$ and the  Lebesgue measure  $|\Omega|<\infty$ that 
		$\frac{1}{m_{\theta,\mu}}\widehat{V_{\theta,\mu}f}\in L^1(\Omega)$.  Then the right-hand side of \eqref{eq:inv}
		defines a continuous function. Consequently, 
		\eqref{eq:inv} holds for all $\mathbf{x}\in \mathbb{R}^2.$
		
		Next we investigate how to represent $\widehat{V_{\theta,\mu}f}$
		by the samples of $V_{\theta,\mu}f$. Applying Lemma \ref{lem:Shannon} to $V_{\theta,\mu}f$
		yields
		\begin{align}	
			\label{AVLTbiaoshi}
			(V_{\theta,\mu}f)(\mathbf{v})=\sum_{\mathbf{k}\in\Z^2}[V_{\theta,\mu}f](\lambda_{\boldsymbol{k}})e^{\mathrm{i}\mathbf{c}\cdot(\mathbf{v}-\lambda_{\boldsymbol{k}})} \mathcal S\left(\boldsymbol{v}\oslash \boldsymbol{h}-\mathbf{k}\right)
		\end{align}
		for all $\mathbf{v}\in\mathbb R^2$,
		where $\mathbf{c}=(c_1,c_2)$ and $\boldsymbol{h}$ are  as in \eqref{canshu123}.
		Consequently, 
		\begin{equation}\label{eq:Vhat}
			\widehat{V_{\theta,\mu}f}(\boldsymbol{\xi})
			=
			h_1h_2
			\chi_{\Omega}(\boldsymbol{\xi})
			\sum_{\mathbf{k}\in\Z^2}
			[V_{\theta,\mu}f](\lambda_{\boldsymbol{k}})e^{-\mathrm{i}\lambda_{\boldsymbol{k}}
				\cdot\boldsymbol{\xi}},
			\quad \boldsymbol{\xi}\in\mathbb{R}^2.
		\end{equation}
		Substituting \eqref{eq:Vhat} into \eqref{eq:inv} and interchanging
		sum and integral, we obtain  \eqref{caiyangdinglidaixian}.
	\end{proof}

	\subsection{Characterization of sampling series in Theorem \ref{thm:sampling-reconstruction}}

	\subsubsection{Riesz basis property of reconstruction system \texorpdfstring
		{$\{\Phi\left(\boldsymbol{x}-\lambda_{\boldsymbol{k}}\right)\}_{\boldsymbol{k}\in \mathbb{Z}^2}$}
		{} and optimal Riesz  bounds}
	This section starts with the following lemma. 
	\begin{lem}\label{lem:min-m}
		As previously,  $\Omega=[a_1,b_1]\times [a_2,b_2].$
		For $j=1,2$, define
		\begin{equation}\label{def:ho,R}
			\rho_j:=
			\begin{cases}
				0, & a_j\le 0\le b_j,\\
				\min\{|a_j|,|b_j|\}, & 0\notin [a_j,b_j],
			\end{cases}
			\ \  \hbox{and} \ \ 
			R_j:=\max\{|a_j|,|b_j|\}.
		\end{equation}
		Additionally, define 
		$ 
		U_0:=\rho_1^2\sin^2\theta,
		U_1:=R_1^2\sin^2\theta,
		V_0:=\rho_2^2\cos^2\theta $ and $
		V_1:=R_2^2\cos^2\theta.
		$
		Then   $\min_{\boldsymbol{\xi}\in \Omega}|m_{\theta,\mu}(\boldsymbol{\xi})|^2$ on $\Omega$ is attained at one of the four points
		\[
		(U_0,V_0),\quad (U_0,V_1),\quad (U_1,V_0),\quad (U_1,V_1), 
		\]
		while  $\max_{\boldsymbol{\xi}\in \Omega}|m_{\theta,\mu}(\boldsymbol{\xi})|^2$  is attained on the boundary of the rectangle $[U_0,U_1]\times[V_0,V_1]$.
	\end{lem}
	
	\begin{proof}
		See Appendix \ref{app:riesz-details} for the proof.
	\end{proof}

	We equip $L^2(\Omega)$ with the normalized inner product
	\begin{align}\label{neijikongjian}
		\langle F,G\rangle_{L^2(\Omega)}
		:=
		\frac{h_1h_2}{(2\pi)^2}
		\int_\Omega F(\mathbf{\xi})\overline{G(\mathbf{\xi})}d\mathbf{\xi},
		\qquad F,G\in L^2(\Omega).
	\end{align}
	Now it is ready to investigate the Riesz basis 
	property of  $\Phi$.

	\begin{theo}\label{th:Riesz-bounds}
		Let $\Phi$ and  $\lambda_{\boldsymbol{k}}$   be as  in  \eqref{TTCX} and  \eqref{caiyangsanshu}, respectively. Then 
		the lattice translates $\left\{\Phi(\cdot-\lambda_{\boldsymbol{k}})\right\}_{\boldsymbol{k}\in\mathbb Z^2}$ form a Riesz basis of $PW_\Omega$. Moreover, the optimal Riesz lower
		and upper bounds are $A_{\mathrm{opt}}:=\frac{h_1h_2}{\displaystyle\max_{\xi\in\Omega}|m_{\theta,\mu}(\xi)|^2}$ and $B_{\mathrm{opt}}:=\frac{h_1h_2}{\displaystyle\min_{\xi\in\Omega}|m_{\theta,\mu}(\xi)|^2}$, respectively. 
	\end{theo}
	
	\begin{proof}
		The proof is  provided in Appendix
		\ref{proof-riesz-basis}.
	\end{proof}
	
	From Theorem \ref{th:Riesz-bounds}, the condition number 
	$\kappa=\sqrt{\frac{B_{\mathrm{opt}}}{A_{\mathrm{opt}}}}=\sqrt{\frac{\max_{\xi\in\Omega}|m_{\theta,\mu}(\xi)|^2}{\min_{\xi\in\Omega}|m_{\theta,\mu}(\xi)|^2}}.$ As shown in the following example,
	for $\Omega$ centered at the origin, $\kappa$
	can be expressed  explicitly.

	\begin{exam}\label{exampleforzuiyoujie}
		Let   $\Omega=[-1,1]^2$. Since $0\in\Omega$,
		we have $\max_{\xi\in\Omega}|m_{\theta,\mu}(\xi)|^2=4/\mu^2$.
		By Lemma  \ref{lem:min-m},
		\[
		\begin{aligned}
			&\min_{\xi\in\Omega}|m_{\theta,\mu}(\xi)|^2\\
			&=
			\min\left\{
			\frac{4\mu^2}{(\mu^2+\sin^2\theta)^2},
			\frac{4}{\mu^2+\cos^2\theta},
			\frac{4(\mu^2+\cos^2\theta)}
			{(\mu^2+(\sin\theta+\cos\theta)^2)
				(\mu^2+(\cos\theta-\sin\theta)^2)}
			\right\}.
		\end{aligned}
		\]
		Hence
		\[
		\begin{aligned}
			&\kappa
			=
			\sqrt{\frac{\max_{\xi\in\Omega}|m_{\theta,\mu}(\xi)|^2}{\min_{\xi\in\Omega}|m_{\theta,\mu}(\xi)|^2}}\\
			&=
			\max\left\{
			1+\frac{\sin^2\theta}{\mu^2},
			\sqrt{1+\frac{\cos^2\theta}{\mu^2}},
			\sqrt{
				\frac{(\mu^2+(\sin\theta+\cos\theta)^2)(\mu^2+(\cos\theta-\sin\theta)^2)}
				{\mu^2(\mu^2+\cos^2\theta)}
			}
			\right\}.
		\end{aligned}
		\]
		Consequently, for each fixed $\theta\in(0,\pi/2)$, $\kappa=1+O(\mu^{-2}),\mu\to\infty$,
		and 
		$\mu^2\kappa=\sin^2\theta+o(1), \ \mu\to0^+$.
		This implies that a larger $\mu$ leads to stronger attenuation and greater system stability of $\left\{\Phi(\cdot-\lambda_{\boldsymbol{k}})\right\}_{\boldsymbol{k}\in\mathbb Z^2}$. 
	\end{exam}
	

	\subsubsection{Adaptive  truncation, uniform convergence and the optimal convergence rate}
	For the sampling series in Theorem \ref{thm:sampling-reconstruction}, a natural problem is to establish its truncation approximation. We consider the adaptive truncation. Specifically, for any fixed $R>0$  and  $\mathbf{x}\in \mathbb{R}^2$, the adaptive truncation is defined as 
	\begin{align}\label{adaptivetruncation}
		S_Rf(\mathbf x)
		:=
		\sum_{\substack{\mathbf k\in\mathbb Z^2\\
				\|\mathbf x-\lambda_{\mathbf k}\|_2\le R}}
		[V_{\theta,\mu}f](\lambda_{\mathbf k})
		\Phi(\mathbf x-\lambda_{\mathbf k}).
	\end{align}
	We   need  to establish an  upper bound for 
	$\sup_{\mathbf{x}\in\mathbb{R}^2}\bigl|f(\mathbf{x})-S_Rf(\mathbf{x})\bigr|$.

	\begin{lem}\label{lem:Phi-decay}
		Set
		\[
		a(\boldsymbol\xi):=\frac{1}{m_{\theta,\mu}(\boldsymbol\xi)},\quad 
		\boldsymbol\xi=(\xi_{1}, \xi_{2})\in\Omega,
		\]
		where $m_{\theta,\mu}(\boldsymbol\xi)$ is as in \eqref{eq:m-theta-mu}. 
		Define
		\begin{align}\label{constants}
			M_0^\Omega:=\|a\|_{L^\infty(\Omega)}, \ \ 
			M_1^\Omega:=\max_{j=1,2}\|\partial_{\xi_j}a\|_{L^\infty(\Omega)},\ \ M_2^\Omega:=\|\partial_{\xi_1}\partial_{\xi_2}a\|_{L^\infty(\Omega)} .
		\end{align}
		Then there exists a constant $C_\Omega>0$
		(depending only on  
		$\Omega$) such that
		\begin{align}\label{shengchengyuanguji}
			|\Phi(x_1,x_2)|
			\le
			C_\Omega
			\frac{M_0^\Omega+M_1^\Omega+M_2^\Omega}
			{(1+|x_1|)(1+|x_2|)}
		\end{align}
		for every $(x_1,x_2)\in \mathbb R^2$. 
	\end{lem}
	
	\begin{proof}
		The proof is provided  in Appendix \ref{app:truncation-details}.
	\end{proof}

	\begin{lem}\label{lem:sampling}
		Let $f\in PW_\Omega$ and  $\lambda_{\boldsymbol{k}}$ be  as in \eqref{caiyangsanshu}. Then
		\begin{align}\label{1111}
			\sum_{\mathbf k\in\mathbb Z^2}
			\left|[V_{\theta,\mu}f](\lambda_{\mathbf k})\right|^2
			\le
			\frac{4}{\mu^2 h_1h_2}
			\|f\|_{L^2(\mathbb R^2)}^2.
		\end{align}
	\end{lem}
	
	\begin{proof}
		Using the  inverse Fourier transform and  \eqref{eq:Fourier identity}, for any  $\mathbf{k}\in\Z^2$ it follows that
		\[
		[V_{\theta,\mu}f](\lambda_{\mathbf k})
		=
		\frac{1}{(2\pi)^2}
		\int_\Omega
		m_{\theta,\mu}(\boldsymbol\xi)\widehat f(\boldsymbol\xi)
		e^{\mathrm{i}\lambda_{\mathbf k}\cdot\boldsymbol\xi}
		d\boldsymbol\xi .
		\]
		Let  $L^2(\Omega)$ be as  in \eqref{neijikongjian}. 
		Then $\{e_{\mathbf{k}}(\boldsymbol{\xi}):=e^{\mathrm{i}\lambda_{\mathbf k}\cdot\boldsymbol{\xi}}\}_{\mathbf{k}\in\Z^2}$ forms an orthonormal basis for 
		$L^2(\Omega)$. 
		By Parseval's identity,
		\[
		\sum_{\mathbf k\in\mathbb Z^2}
		\left|[V_{\theta,\mu}f](\lambda_{\mathbf k})\right|^2
		=
		\frac{1}{(h_1h_2)^2}
		\|m_{\theta,\mu}\widehat f\|_{\Omega}^2
		=
		\frac{1}{h_1h_2}
		\frac{1}{(2\pi)^2}
		\int_\Omega
		|m_{\theta,\mu}(\boldsymbol\xi)\widehat f(\boldsymbol\xi)|^2
		d\boldsymbol\xi.
		\]
		Using the bound $|m_{\theta,\mu}(\boldsymbol{\xi})|\le \|m_{\theta,\mu}\|_\infty$  together with  Plancherel's theorem  yields the claim.
	\end{proof}
	
	The following lemma will be  helpful for establishing the optimal convergence rate in Theorem \ref{thm:truncation}.

	\begin{lem}\label{lem:Phi-lower-bound}
		Let $\mathbf x^0=\left(\frac{h_1}{2},0\right)$ and $t_n=\frac{h_1}{2}-h_1n$ with $ n\in\mathbb N$.
		Then there exists $N_0\in\mathbb N$ such that $|\Phi(t_n,0)|\gtrsim 1/n$ whenever $n\ge N_0$. 
	\end{lem}
	
	\begin{proof}
		See Appendix \ref{proof-fuzhu-zuiyou}.
	\end{proof}
	
	\begin{theo}\label{thm:truncation}
		For any $f\in PW_{\Omega}$ and  $R>0$,
		$S_Rf(\mathbf x)$ is defined in \eqref{adaptivetruncation}.
		Then 
		$\sup_{\mathbf{x}\in\mathbb{R}^2}\bigl|f(\mathbf{x})-S_Rf(\mathbf{x})\bigr|\lesssim R^{-1/2}$ as $R\rightarrow \infty$. Moreover, the convergence rate $1/2$ is optimal, namely, there exists $f\in PW_\Omega$ such that, for every $\delta>0$, $\sup_{\mathbf x\in\mathbb R^2}
		|f(\mathbf x)-S_Rf(\mathbf x)|\ne O(R^{-1/2-\delta})$ as $R\rightarrow \infty$.
	\end{theo}
	
	\begin{proof}
		The complete argument is given in Appendix \ref{proof-zuiyou}.
	\end{proof}
	
	\begingroup
	\raggedbottom
	\subsection{Numerical simulation}\label{sec:numerical}
	The purpose of this section is to verify the effectiveness of the
	approximation formula in \eqref{adaptivetruncation}.
	Without loss of generality, we take $\Omega=[-1,1]^2$,
	$\theta=\pi/4$ and $\mu\in\{0.08, 0.1, 0.12\}$.
	As an example, we consider
	\begin{align}\label{xianhaobiaodashi}
		f(\boldsymbol{x})
		=
		\sum_{\boldsymbol{\lambda}\in\mathcal K}
		c_{\boldsymbol{\lambda}}
		\mathcal S\left(\frac{\boldsymbol{x}}{\pi}-\boldsymbol{\lambda}\right),
	\end{align}
	where $\mathcal S(\mathbf{x}):=\operatorname{sinc}(x_1)\operatorname{sinc}(x_2)$  with $\boldsymbol{x}=(x_1,x_2)$, and $ \boldsymbol{\lambda}=(\lambda_1,\lambda_2)\in\mathcal K:=\{0,\ldots,20\}^2$.
	The graph of $f$ on $[0, 20\pi]^{2}$ is plotted in 
	the first column of Fig.\ref{fig:bandlimited-noiseless}.
	Following  \eqref{adaptivetruncation}, we approximate    $f$ by 
	\begin{align}\label{shiyanbufenchongjiangongshi}
		S_Rf(\mathbf{x})=
		\sum_{\substack{\boldsymbol{\gamma}\in\pi\mathbb Z^2\\
				\|\mathbf{x}-\bm{\gamma}\|_2\le R}}
		V_{\theta,\mu}f(\bm{\gamma})\Phi(\mathbf{x}-\bm{\gamma}),
	\end{align}
	where $R/\pi\in\{15,30,200,1000,10000\}$. 
	Since  $c_{\boldsymbol{\lambda}}=f(\gamma_1,\gamma_2)$, 
	we assess the reconstruction effectiveness  using the relative reconstruction error
	$$
	\mathrm{RES}:=\frac{\|C-\widetilde C\|_F}{\|C\|_F},
	$$
	where $\widetilde C=(S_Rf(\gamma_1,\gamma_2))_{\gamma_1,\gamma_2=0}^{20\pi}$
	and $\|\cdot\|_F$ denotes the Frobenius norm. 
	Correspondingly, we need to reconstruct $f$  at every  $\bm{x}\in D:=\{0, \pi, \ldots, 20\pi\}^2$.
	Following  \eqref{eq:Fourier identity}, we have 
	\begin{align}\label{FFTjisuanAVL}
		V_{\theta,\mu}f(\boldsymbol{x})
		=
		\frac{1}{(2\pi)^2}
		\int_{[-1,1]^2}
		m_{\theta,\mu}(\boldsymbol{\xi})
		\widehat{f}(\boldsymbol{\xi})
		e^{i\boldsymbol{x}\cdot\boldsymbol{\xi}}
		d\boldsymbol{\xi}.
	\end{align}
	For this we use the  2D fast Fourier transform (FFT) to compute 
	$\{V_{\theta,\mu}f(\bm{\gamma})\}$.

	\subsubsection{Noiseless setting}
	For each fixed $\mu$, Table \ref{tab:bandlimited-truncation}
	suggests that RES decays at the rate $O(R^{-1})$ as $R$ increases.
	This does not contradict Theorem \ref{thm:truncation}, since
	the rate $O(R^{-1/2})$ established therein is optimal 
	(sharp) over
	the entire space $PW_{\Omega}$, whereas some  individual functions
	may exhibit faster decay.
	On the other hand, for a fixed $R>0$, the reconstruction
	accuracy improves as $\mu$ increases.
	This can be explained as follows.
	Example \ref{exampleforzuiyoujie} shows that the optimal
	Riesz upper bound $B_{\mathrm{opt}}$ for
	$\{\Phi(\cdot-\bm{\gamma})\}_{\bm{\gamma}\in\pi\mathbb Z^2}$
	decreases as $\mu$ increases.
	The truncation  (or  omitting some terms
	of $\{\Phi(\cdot-\bm{\gamma})\}_{\bm{\gamma}\in\pi\mathbb Z^2}$) are independent
	of $\mu$, and \eqref{Rieszbudengshi} yields a smaller upper bound
	on the truncation error as $\mu$ increases.
	The graphs of $S_Rf$ on $[0,20\pi]^2$ for different values of
	$\mu$ and $R$ are shown in panels $(a_2)$--$(a_4)$,
	$(b_2)$--$(b_4)$, and $(c_2)$--$(c_4)$ of
	Fig.\ref{fig:bandlimited-noiseless}.
	\begin{table}[!htbp]
		\centering
		\caption{RES for noiseless recovery.}
		\label{tab:bandlimited-truncation}
		\renewcommand{\arraystretch}{1.10}
		\begin{tabular*}{0.85\textwidth}{@{\extracolsep{\fill}}cccc@{}}
			\specialrule{\heavyrulewidth}{0pt}{0pt}    \diagbox[width=4.6em,height=2.3em,linewidth=0.45pt,trim=l,innerleftsep=3pt,innerrightsep=3pt]{$R/\pi$}{$\mu$} & $0.08$ & $0.1$ & $0.12$ \\
			\specialrule{\lightrulewidth}{0pt}{\belowrulesep}
			$15$ & 6.1225e-2 & 5.8156e-2 & 5.5602e-2 \\
			$30$ & 2.3849e-2 & 2.2868e-2 & 2.1927e-2 \\
			$200$ & 3.0970e-3 & 2.9650e-3 & 2.8391e-3 \\
			$1000$ & 5.6728e-4 & 5.4298e-4 & 5.1981e-4 \\
			$10000$ & 4.8021e-5 & 4.5961e-5 & 4.3998e-5 \\
			\bottomrule
		\end{tabular*}
	\end{table}
	\begin{figure}[!tbp]
		\centering
		\includegraphics[width=\textwidth]{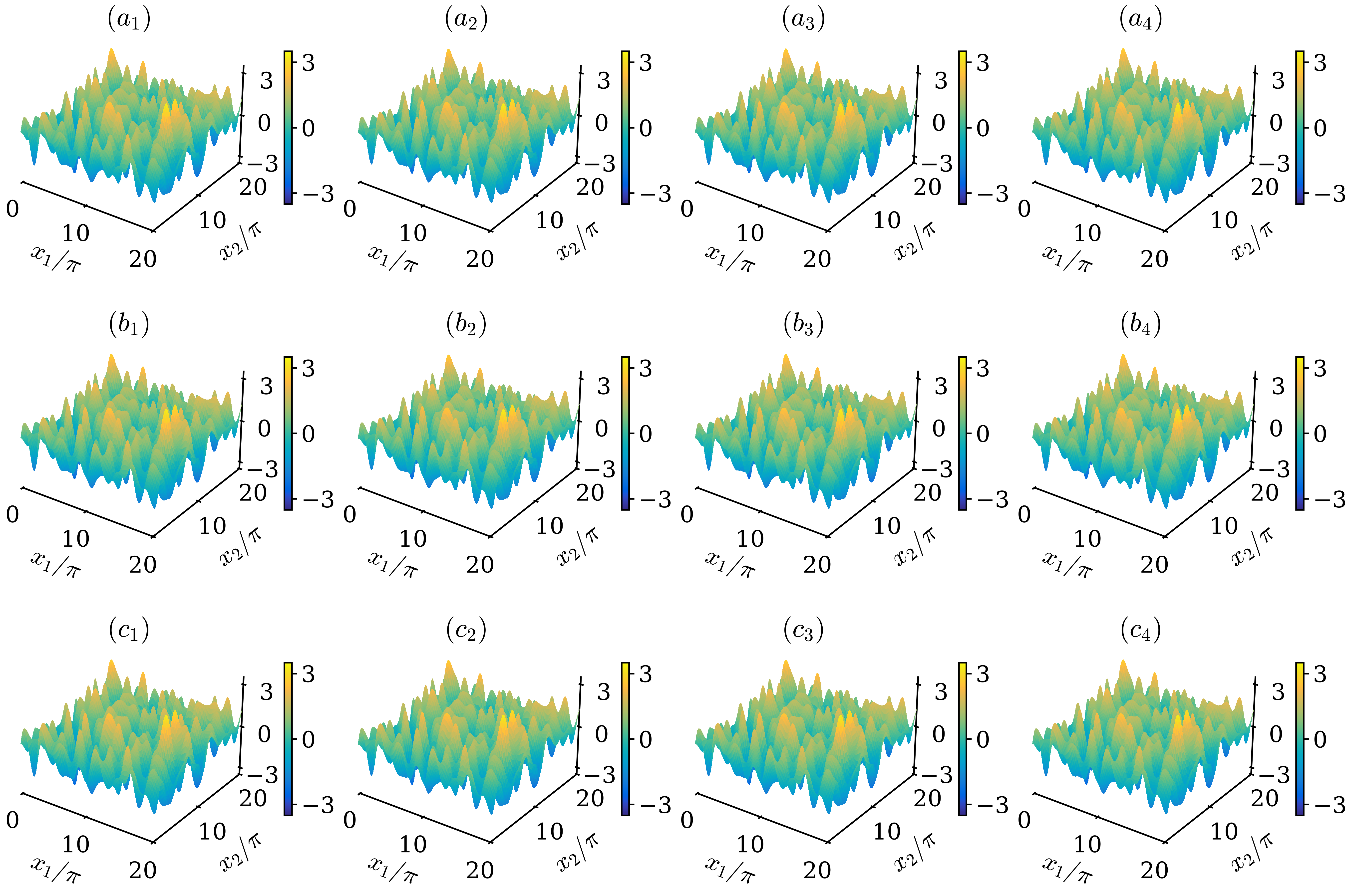}
		\caption{Function reconstructions in the noiseless setting  for $\theta=\pi/4$. The first column shows
			the original version; columns 2--4 show reconstructions versions for $R=15\pi,30\pi,200\pi$. Rows
			$(a_1)$--$(a_4)$, $(b_1)$--$(b_4)$, and $(c_1)$--$(c_4)$
			correspond to $\mu=0.08,0.1,0.12$, respectively. }
		\label{fig:bandlimited-noiseless}
	\end{figure}

	\subsubsection{Noisy setting}\label{nosysetting}
	The purpose of this subsection is to check the   robustness of \eqref{shiyanbufenchongjiangongshi} to the  additive Gaussian noises. For simplicity, we just address the case of $R=10000\pi$. The sample set for all $\bm{x}\in D$
	is 
	$\{V_{\theta,\mu}f(\bm{\gamma)}\}_{\boldsymbol{\gamma}\in\mathcal A_{R}}$
	where 
	$$
	\mathcal A_{R}:=\left\{\bm{\gamma}\in\pi\mathbb Z^2:
	\inf_{\mathbf{x}\in D}\|\bm{\gamma}-\mathbf{x}\|_2\le R\right\}.
	$$
	

	Suppose that the samples $V_{\theta,\mu}f(\boldsymbol{\gamma})$,
	$\boldsymbol{\gamma}\in\mathcal A_{R}$, are corrupted by the independent
	and identically distributed (i.i.d.) Gaussian noise
	$\epsilon_{\boldsymbol{\gamma}}\sim N(0,\sigma^2)$.
	The observed samples are therefore
	$V_{\theta,\mu}f(\boldsymbol{\gamma})
	+\epsilon_{\boldsymbol{\gamma}}$, where
	the variance $\sigma^2$ is determined by the
	signal-to-noise ratio:
	\[
	\operatorname{SNR}\;(\mathrm{dB})
	=10\log_{10}\left(
	\frac{\sum_{\boldsymbol{\gamma}\in\mathcal A_{R}}
		|V_{\theta,\mu}f(\boldsymbol{\gamma})|^2}
	{\#\mathcal A_{R}\sigma^2}
	\right).
	\]
	We perform $100$  trials  for each $\mu\in\{0.08,0.1,0.12\}$
	and report the mean RESs  in Table \ref{tab:bandlimited-noise}.
	For fixed values of SNR and $R$, the RES decreases as $\mu$ increases.
	This observation is consistent with the condition numbers given  in
	Example \ref{exampleforzuiyoujie}, which satisfy
	\[
	\left.\kappa\right|_{\mu=0.08}
	\approx79.125>\left.\kappa\right|_{\mu=0.1}
	\approx51>\left.\kappa\right|_{\mu=0.12}
	\approx35.72.
	\]
	Without bias, for  SNR=30 and  each pair ($\mu$, $R$) we show 
	the  reconstruction version  of $f$
	on $[0, 20\pi]^2$   in 
	Fig.\ref{fig:bandlimited-noise-images}.

	\begin{table}[!htbp]
		\centering
		\caption{Mean RES of results of   100 trials at $R=10000\pi$.}
		\label{tab:bandlimited-noise}
		\renewcommand{\arraystretch}{1.10}
		\begin{tabular*}{0.85\textwidth}{@{\extracolsep{\fill}}cccc@{}}
			\specialrule{\heavyrulewidth}{0pt}{0pt}
			\diagbox[width=8.5em,height=3em,linewidth=0.45pt,trim=l,innerleftsep=15pt,innerrightsep=3pt]{SNR (dB)}{\raisebox{-0.85em}[\height][\depth]{$\mu$}}
			& $0.08$ & $0.1$ & $0.12$ \\
			\specialrule{\lightrulewidth}{0pt}{\belowrulesep}
			$15$ & 8.1342e-4 & 6.3834e-4 & 5.4687e-4 \\
			$20$ & 4.5915e-4 & 3.6101e-4 & 3.0968e-4 \\
			$25$ & 2.6125e-4 & 2.0658e-4 & 1.7792e-4 \\
			$30$ & 1.5221e-4 & 1.2228e-4 & 1.0648e-4 \\
			\bottomrule
		\end{tabular*}
	\end{table}	
	\begin{figure}[!tbp]
		\centering
		\includegraphics[width=\textwidth]{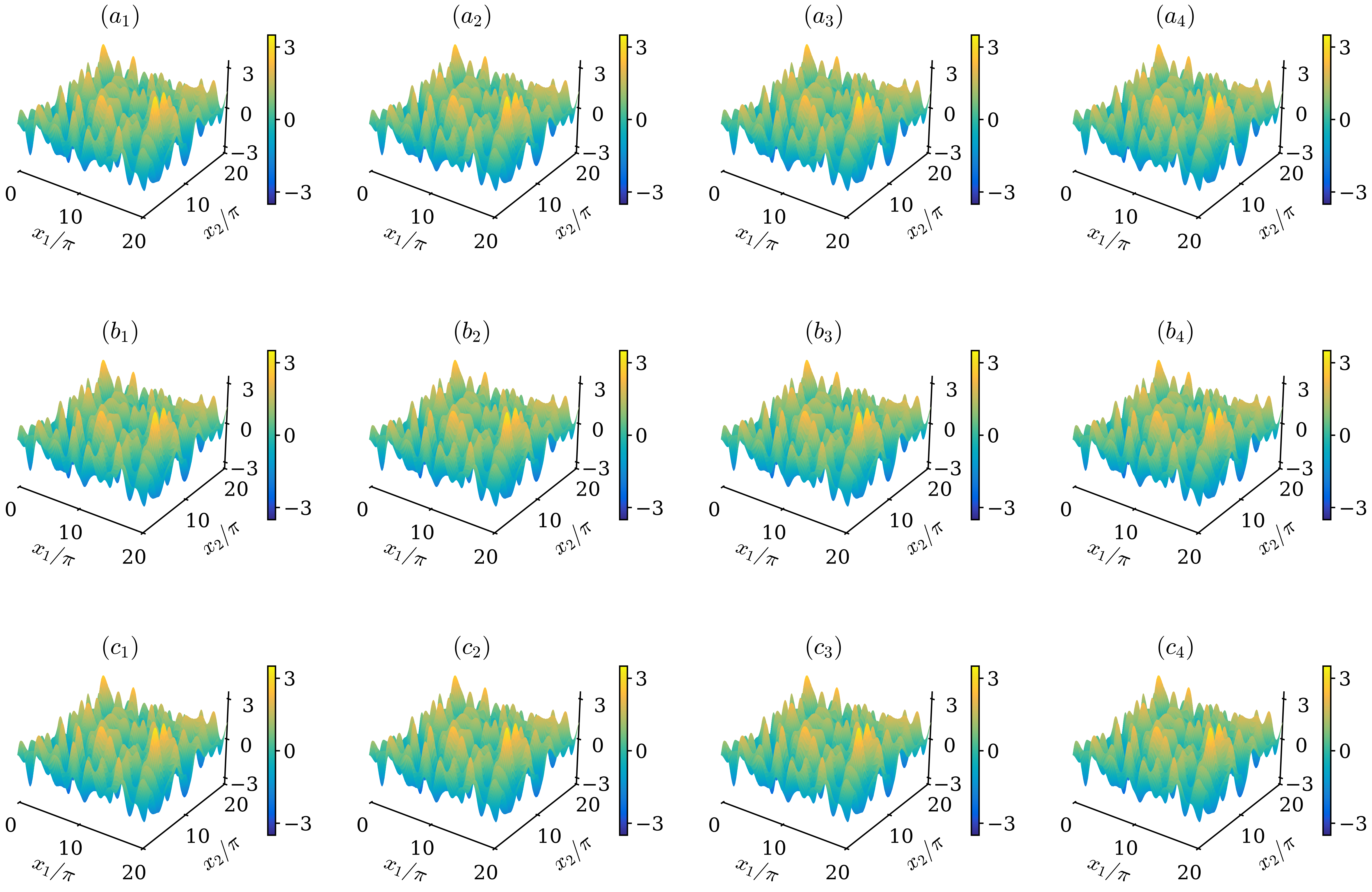}
		\caption{Function recovery from noisy AV-LT samples at
			SNR $30\,\mathrm{dB}$ and $\theta=\pi/4$.
			Panels  $(a_1)$--$(a_4)$, $(b_1)$--$(b_4)$, and
			$(c_1)$--$(c_4)$ correspond to $\mu=0.08,0.1,0.12$, respectively.
			The first column shows the original image; columns 2--4 show reconstructions
			at $R=200\pi,1000\pi,10000\pi$, respectively.}
		\label{fig:bandlimited-noise-images}
	\end{figure}
	
	\begin{rem}
		The numerical experiments involve three sources of error: (1) the truncation in \eqref{shiyanbufenchongjiangongshi}; (2) the approximate computation of the data $\{V_{\theta,\mu}f(\bm{\gamma})\}$ using the FFT; (3) noise contamination. The numerical results demonstrate the stability of the reconstruction method in the presence of these three sources of  error.
	\end{rem}

	\subsection{Comparison with the  Tikhonov regularization-based reconstruction}
	By \eqref{xianhaobiaodashi}, we have 
	\begin{align}\label{xianhaobiaodashinew}
		V_{\theta,\mu}f(\mathbf{x})=
		\sum_{\boldsymbol{\lambda}\in\mathcal K}
		c_{\boldsymbol{\lambda}}
		V_{\theta,\mu}\left[
		\mathcal S\left(\frac{\cdot}{\pi}-\boldsymbol \lambda\right)
		\right](\mathbf{x}),
	\end{align}
	where $\mathcal K=\{ 0, 1, \ldots, 20\}^{2}$.
	Define a vector  $\bm{c}=(c_{\boldsymbol{\lambda}})_{\boldsymbol{\lambda}\in \mathcal K}$.
	We  reconstruct 
	$\bm{c}$
	through the regularization equation 
	\begin{align}
		(Y^{*}Y+\delta I)\bm{c}=Y^{*}(V_{\theta,\mu}f(\pi\boldsymbol{\lambda})+\epsilon_{\mathbf{\lambda}})_{\mathbf{\lambda}\in \mathcal K}
	\end{align}
	where $Y:=\Big(V_{\theta,\mu}\left[
	\mathcal S\left(\frac{\cdot}{\pi}- \boldsymbol \lambda \right)
	\right](\pi\mathbf{k})\Big)_{\mathbf{\lambda},\mathbf{k}\in \mathcal K}$ is a $441\times 441$ matrix
    with $Y^{*}$ its conjugate transpose, 
	$\epsilon_{\mathbf{\lambda}}$
	is the additive Gaussian noise defined as in section \ref{nosysetting}, and $\delta$ is regularization parameter  selected  through the
	Morozov discrepancy principle. For each  pair $(\mu, \hbox{SNR})\in \{0.08, 0.1, 0.12\} \times \{15, 20, 25, 30\}$, the simulation is conducted in the noisy  setting 
	for 100 trials. Their mean RESs are reported in Table 
	\ref{tab:sinc21-regularized}. 
	Without bias, for each pair ($\mu$, SNR) we show 
	a reconstruction version  of $f$
	on $[0, 20\pi]^2$   in 
	Fig.\ref{fig:sinc21-regularized-surfaces}.
	Comparing Tables  \ref{tab:bandlimited-noise} and \ref{tab:sinc21-regularized}, Figs.\ref{fig:bandlimited-noise-images} and \ref{fig:sinc21-regularized-surfaces}, we 
	found that the direct reconstruction method \eqref{shiyanbufenchongjiangongshi} outperforms the regularization method. 
	
	\begin{figure}[!tbp]
		\centering
		\includegraphics[width=\textwidth]{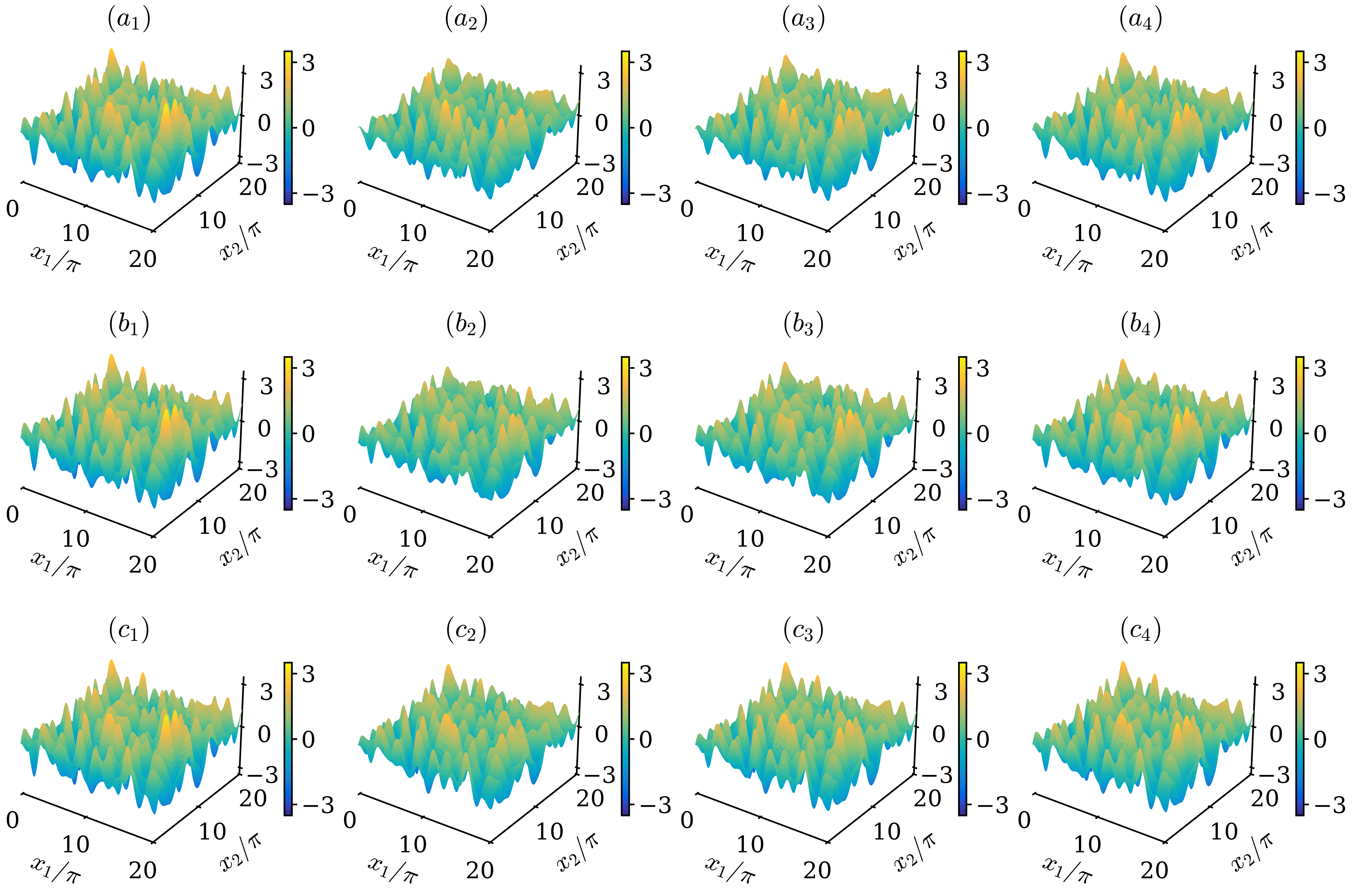}
		\caption{Function reconstruction by Tikhonov regularization
			from noisy AV-LT samples at $\theta=\pi/4$.
			Rows $(a_1)$--$(a_4)$, $(b_1)$--$(b_4)$, and
			$(c_1)$--$(c_4)$ correspond to $\mu=0.08,0.1,0.12$, respectively.
			The first column shows $f$; columns 2--4 show reconstructions
			at SNRs $15$, $20$, and $25$ dB, respectively.}
		\label{fig:sinc21-regularized-surfaces}
	\end{figure}
	
	\begin{table}[!htbp]
		\centering
		\caption{Mean RES of  Tikhonov regularization-based  reconstruction.}
		\label{tab:sinc21-regularized}
		\renewcommand{\arraystretch}{1.10}
		\begin{tabular*}{0.85\textwidth}{@{\extracolsep{\fill}}cccc@{}}
			\specialrule{\heavyrulewidth}{0pt}{0pt}
			\diagbox[width=8.5em,height=3em,linewidth=0.45pt,trim=l,innerleftsep=15pt,innerrightsep=3pt]{SNR (dB)}{\raisebox{-0.85em}[\height][\depth]{$\mu$}}
			& $0.08$ & $0.1$ & $0.12$ \\
			\specialrule{\lightrulewidth}{0pt}{\belowrulesep}
			$15$ & 4.5551e-1 & 4.1503e-1 & 3.8785e-1 \\
			$20$ & 3.1708e-1 & 2.8109e-1 & 2.5797e-1 \\
			$25$ & 2.0974e-1 & 1.8089e-1 & 1.6330e-1 \\
			$30$ & 1.3247e-1 & 1.1157e-1 & 9.9434e-2 \\
			\bottomrule
		\end{tabular*}
	\end{table}

	\par\endgroup

	\section{An approximate reconstruction formula for Sobolev spaces}\label{Sobolev space jin si}
	Recall that the  sampling  in Theorem \ref{thm:sampling-reconstruction} holds for band-limited functions. In this section, we aim to establish an approximate   reconstruction formula  for functions in  
	$H^{\lambda}(\mathbb{R}^2)$ from their AV-LT samples, where $\lambda>2$. 
	

	\subsection{Answer to \textbf{Q5}: a sampling-based approximation scheme}
	For $f\in H^\lambda(\mathbb{R}^2)$  with $\lambda>2$,  we choose as a bandlimited proxy the $L^2$-orthogonal projection
	$P_N f\in PW_{\Omega_N}$   defined by 
	\begin{align}\label{daixianbijin}
		\widehat{P_N f}(\boldsymbol{\xi})=\chi_{\Omega_N}(\boldsymbol{\xi})\widehat f(\boldsymbol{\xi}),\ \ \boldsymbol{\xi}\in\mathbb{R}^2,   
	\end{align}
	where  $\Omega_N=[-2^N\pi,2^N\pi]^2$,     $N$  denotes the scale level, and $\chi_{\Omega_N}$ is the indicator function of $\Omega_N$.
	
	Since $P_N f\in PW_{\Omega_N}$, \eqref{caiyangdinglidaixian}      yields 
	its  reconstruction from the  AV-LT samples. Specifically,  
	\begin{align}\label{daixinabijin1}
		(P_N f)(\mathbf{x})
		=\sum_{\mathbf{k}\in\Z^2}[V_{\theta,\mu}(P_N f)]\!\left(2^{-N}\mathbf{k}\right)
		\Phi\bigl(\mathbf{x}-2^{-N}\mathbf{k}\bigr),\ \ \mathbf{x}\in\mathbb{R}^2.  
	\end{align}
	In practice, however, the available measurements are  
	$\{[V_{\theta,\mu}f](2^{-N}\mathbf{k})\}_{\mathbf{k}\in\Z^2}$ instead of 
	$\{[V_{\theta,\mu}(P_N f)](2^{-N}\mathbf{k})\}_{\mathbf{k}\in\Z^2}$.
	This, together with \eqref{daixinabijin1} and \eqref{daixianbijin},  motivates the following approximation scheme 
	\begin{align}\label{sobolevseries}
		f\approx f_N^\#(\mathbf{x})
		:=\sum_{\mathbf{k}\in\Z^2}[V_{\theta,\mu}f](2^{-N}\mathbf{k})
		\Phi\bigl(\mathbf{x}-2^{-N}\mathbf{k}\bigr),\ \  \mathbf{x}\in\mathbb{R}^2.
	\end{align}
	
	\subsection{Approximation error}
	The aim of this subsection is to establish an error estimate for $\|f-f_N^\#\|_{2}$. For every $f\in H^\lambda(\mathbb{R}^2)$ and $N\in\N$,
	\begin{align}\label{wuchafenjie}
		f-f_N^\#=(I-P_N)f+\bigl(P_N f-f_N^\#\bigr),
	\end{align}
	where 
	\begin{equation}\label{eq:alg-error}
		P_N f-f_N^\#
		=-\sum_{\mathbf{k}\in\Z^2}[V_{\theta,\mu}((I-P_N)f)]\!\left(2^{-N}\mathbf{k}\right)
		\Phi\bigl(\mathbf{x}-2^{-N}\mathbf{k}\bigr).
	\end{equation}
	Then we just need to estimate $\|(I-P_N)f\|_{2}$
	and $\|P_N f-f_N^\#\|_{2}$.

	\begin{lem}\label{lem:inv-m-growth}
		There exists a constant $C_{\mu}>0$ such that for all $\boldsymbol{\xi}\in\Omega_N$,
		\[
		\sup_{\boldsymbol{\xi}\in\Omega_N}\Bigl|\frac{1}{m_{\theta,\mu}(\boldsymbol{\xi})}\Bigr|
		\le C_{\mu}\bigl(1+(2^N\pi)^2\bigr).
		\]
	\end{lem}
	
	\begin{proof}
		By \eqref{eq:m-theta-mu}, we have
		$$
		m_{\theta,\mu}(\boldsymbol{\xi})
		=
		\frac{2(\mu-\mathrm{i}\xi_2\cos\theta)}
		{(\mu-\mathrm{i}\mathbf d_+(\theta)\cdot\boldsymbol{\xi})
			(\mu-\mathrm{i}\mathbf d_-(\theta)\cdot\boldsymbol{\xi})}.
		$$
		Since $\|\mathbf d_\pm(\theta)\|_2=1$ and
		$|\mu-\mathrm{i}\xi_2\cos\theta|\ge\mu>0$, it follows that
		$$
		\left|\frac{1}{m_{\theta,\mu}(\boldsymbol{\xi})}\right|
		=
		\frac{
			|\mu-\mathrm{i}\mathbf d_+(\theta)\cdot\boldsymbol{\xi}|
			|\mu-\mathrm{i}\mathbf d_-(\theta)\cdot\boldsymbol{\xi}|
		}
		{
			2|\mu-\mathrm{i}\xi_2\cos\theta|
		}
		\le
		\frac{\mu^2+\|\boldsymbol{\xi}\|_2^2}{2\mu}.
		$$
		For $\boldsymbol{\xi}\in\Omega_N$, we have
		$\|\boldsymbol{\xi}\|_2^2\le 2(2^N\pi)^2$. Therefore,
		$$
		\sup_{\boldsymbol{\xi}\in\Omega_N}
		\left|\frac{1}{m_{\theta,\mu}(\boldsymbol{\xi})}\right|
		\le
		\frac{\mu}{2}+\frac{(2^N\pi)^2}{\mu}
		\le
		C_\mu\bigl(1+(2^N\pi)^2\bigr),
		$$
		where $C_\mu:=\max\{\mu/2,\mu^{-1}\}$.
		This completes the proof.
	\end{proof}


	Note that
	$\operatorname{supp}\bigl((f-P_Nf)^{\wedge}\bigr)
	\subseteq\Omega_N^c$.
	By \eqref{eq:Fourier identity} and
	Remark \ref{rem:m-properties},
	$\operatorname{supp}\bigl([V_{\theta,\mu}(f-P_Nf)]^{\wedge}\bigr)
	\subseteq\Omega_N^c$ as well.   Motivated by this observation, we derive an estimate for the sampling  sequence $\{g(2^{-N}\mathbf{k})\}_{\mathbf{k}\in\mathbb{Z}^2}$ as follows. Such an estimate helps us bound $\|P_N f-f_N^\#\|_{2}$.
	
	\begin{lem}\label{lem:aliasing}
		Let $\lambda>2$. 
		Then there exists a constant  $C_{\lambda}>0$ (being independent of $g$) such that 
		\[
		\sum_{\mathbf{k}\in\Z^2}{\left| g\!\left( 2^{-N}\mathbf{k} \right) \right|^2}
		\le C_\lambda2^{-2(\lambda-1)N}\|g\|_{H^\lambda(\mathbb{R}^2)}^2,
		\]
		for  any $g\in H^\lambda(\mathbb{R}^2)$ satisfying  $\operatorname{supp}(\widehat g)\subset \Omega_N^c$.  
	\end{lem}
	
	\begin{proof}
		For $\boldsymbol{\xi}\in\Omega_N$,  define  the $2^{N+1}\pi\Z^2$-periodization of $\widehat{g}$ as 
		\[
		G_N(\boldsymbol{\xi}):=\sum_{\mathbf{n}\in\Z^2}\widehat g\bigl(\boldsymbol{\xi}+2^{N+1}\pi \mathbf{n}\bigr).
		\]
		Splitting $\mathbb{R}^2$ into translates of $\Omega_N$ and using that
		$e^{\mathrm{i}\frac{\mathbf{k}}{2^N}\cdot 2^{N+1}\pi \mathbf{n}}=e^{\mathrm{i}2\pi \mathbf{k}\cdot \mathbf{n}}=1$, we obtain
		\[
		g\!\left(2^{-N}\mathbf{k}\right)
		=\frac{1}{(2\pi)^2}\int_{\Omega_N}G_N(\boldsymbol{\xi})e^{\mathrm{i}\frac{\mathbf{k}}{2^N}\cdot\boldsymbol{\xi}}d\boldsymbol{\xi}.
		\]
		That is,  $\{g\left(2^{-N}\mathbf{k}\right)\}_{\mathbf{k}\in\Z^2}$ are Fourier coefficients of $G_N$ on $\Omega_N$.
		By Parseval identity, we have 
		\begin{align}\label{HT12}
			\sum_{\mathbf{k}\in\Z^2}{\left| g\!\left( 2^{-N}\mathbf{k} \right) \right|^2}
			=\frac{2^{2N}}{(2\pi)^2}\int_{\Omega_N}|G_N(\boldsymbol{\xi})|^2d\boldsymbol{\xi}.
		\end{align}
		Clearly,  
		$\widehat{g}\equiv0$ on $\Omega_N$.  
		Hence $G_N(\boldsymbol{\xi})=\sum_{\mathbf{n}\neq\mathbf{0}}\widehat g(\boldsymbol{\xi}+2^{N+1}\pi \mathbf{n})$.
		By Cauchy-Schwarz inequality and the fact that
		$\left\|\boldsymbol{\xi}\right\|_2\le \sqrt22^N\pi$ on $\Omega_N$, we obtain
		
		\begin{align}\notag
			|G_N(\boldsymbol{\xi})|^2
			&=\left|\sum_{\mathbf{n}\neq\mathbf{0}}\left[(1+\left\|2^{N+1}\pi \mathbf{n}\right\|_2^2)^{-\frac{\lambda}{2}}\right]\left[(1+\left\|2^{N+1}\pi \mathbf{n}\right\|_2^2)^{\frac{\lambda}{2}}\widehat g(\boldsymbol{\xi}+2^{N+1}\pi \mathbf{n})\right]\right|^2\\
			&\label{FFT123}\le \Bigl(\sum_{\mathbf{n}\neq\mathbf{0}}(1+\left\|2^{N+1}\pi \mathbf{n}\right\|_2^2)^{-\lambda}\Bigr)
			\Bigl(\sum_{\mathbf{n}\neq\mathbf{0}}(1+\left\|\boldsymbol{\xi}+2^{N+1}\pi \mathbf{n}\right\|_2^2)^\lambda|\widehat g(\boldsymbol{\xi}+2^{N+1}\pi \mathbf{n})|^2\Bigr). 
		\end{align}
		For $\lambda>2$, we have  $\sum_{\mathbf{n}\neq\mathbf{0}}\left\|\mathbf{n}\right\|_2^{-2\lambda}<\infty$, from which there exists $C_\lambda>0$ such that $\sum_{\mathbf{n}\neq\mathbf{0}}(1+\left\|2^{N+1}\pi \mathbf{n}\right\|_2^2)^{-\lambda}\leq C_\lambda2^{-2\lambda N+2}\pi^2$.
		Combining this and \eqref{FFT123} we have 
		\begin{align}\label{HHT}
			\int_{\Omega_N}|G_N(\boldsymbol{\xi})|^2d\boldsymbol{\xi}
			\le C_\lambda2^{-2\lambda N+2}\pi^2\int_{\mathbb{R}^2}(1+\left\|\boldsymbol{\eta}\right\|^2)^\lambda|\widehat g(\boldsymbol{\eta})|^2d\boldsymbol{\eta}.
		\end{align}
		The proof is completed by \eqref{HT12} and \eqref{HHT}. 
	\end{proof}
	
	\begin{theo}\label{thm:sobolev-rate}
		Let $\lambda>2$. Then  there exists $C>0$ such that for all $f\in H^\lambda(\mathbb{R}^2)$ and all $N\in\N$,
		\begin{align}\label{bijinwucha}
			\|f-f_N^\#\|_{L^2(\mathbb{R}^2)}
			\le C\bigl(1+(2^N\pi)^2\bigr)^{1-\frac{\lambda}{2}}
			\|f\|_{H^\lambda(\mathbb{R}^2)}.
		\end{align}
		In particular, $f_N^\#$
		converges exponentially to
		$f$ in $L^2(\mathbb{R}^2)$ as $N\to\infty$.
	\end{theo}
	
	\begin{proof}
		Obviously, we have $\|f-f_N^\#\|_2\le \|(I-P_N)f\|_2+\|P_N f-f_N^\#\|_2.$
		Since $\left\|\boldsymbol{\xi}\right\|_2\ge 2^N\pi$ on $\mathbb{R}^2\setminus\Omega_N$, then
		\[
		\|(I-P_N)f\|_2^2
		=\frac{1}{(2\pi)^2}\int_{\mathbb{R}^2\setminus\Omega_N}|\widehat f(\boldsymbol{\xi})|^2d\boldsymbol{\xi}
		\le (1+(2^N\pi)^2)^{-\lambda}\|f\|_{H^\lambda(\mathbb{R}^2)}^2.
		\]
		Hence $\|(I-P_N)f\|_2\le (1+(2^N\pi)^2)^{-\lambda/2}\|f\|_{H^\lambda(\mathbb{R}^2)}$.
		
		From \eqref{eq:alg-error},
		\[
		f_N^\#-P_N f= \sum_{\mathbf{k}\in\Z^2}c_{\mathbf{k}}
		\Phi\bigl(\mathbf{x}-2^{-N}\mathbf{k}\bigr), \ \
		c_{\mathbf{k}}:=[V_{\theta,\mu}((I-P_N)f)]\!\left(2^{-N}\mathbf{k}\right).
		\]
		By Theorem \ref{th:Riesz-bounds} and Lemma~\ref{lem:inv-m-growth},
		\[
		\|P_N f-f_N^\#\|_2
		\le 2^{-N}\sup_{\boldsymbol{\xi}\in\Omega_N}\Bigl|\frac{1}{m_{\theta,\mu}(\mathbf\xi)}\Bigr|\|c_{\mathbf{k}}\|_{\ell^2}
		\le C2^{-N}\bigl(1+(2^N\pi)^2\bigr)\|c_{\mathbf{k}}\|_{\ell^2}.
		\]
		Now let $g:=V_{\theta,\mu}(I-P_N)f$. By  \eqref{eq:Fourier identity}, $\widehat g(\boldsymbol{\xi})=m_{\theta,\mu}(\boldsymbol{\xi})(1-\chi_{\Omega_N}(\boldsymbol{\xi}))\widehat f(\boldsymbol{\xi}),$
		so $\operatorname{supp}(\widehat{g})\subseteq\Omega_N^c$. Moreover $|m_{\theta,\mu}|\le 2/\mu$ implies
		$\|g\|_{H^\lambda(\mathbb{R}^2)}\le (2/\mu)\|f\|_{H^\lambda(\mathbb{R}^2)}$. Hence Lemma~\ref{lem:aliasing} gives
		\[
		\|c_{\mathbf{k}}\|_{\ell^2}^2=\sum_{\mathbf{k}\in\Z^2}\left|g\!\left(2^{-N}\mathbf{k}\right)\right|^2
		\le C_\lambda2^{-2(\lambda-1)N}\|g\|_{H^\lambda(\mathbb{R}^2)}^2
		\le C2^{-2(\lambda-1)N}\|f\|_{H^\lambda(\mathbb{R}^2)}^2.
		\]
		
		Therefore,
		\[
		\|P_N f-f_N^\#\|_2
		\le C2^{-N}\bigl(1+(2^N\pi)^2\bigr)2^{-(\lambda-1)N}\|f\|_{H^\lambda(\mathbb{R}^2)}
		\le C\bigl(1+(2^N\pi)^2\bigr)^{1-\frac{\lambda}{2}}\|f\|_{H^\lambda(\mathbb{R}^2)}.
		\]
		This completes the proof.
	\end{proof}
	
	\subsection{Numerical simulation}\label{sec:numerical2}
\begin{figure}[!tbp]
		\centering
		\includegraphics[width=\textwidth]{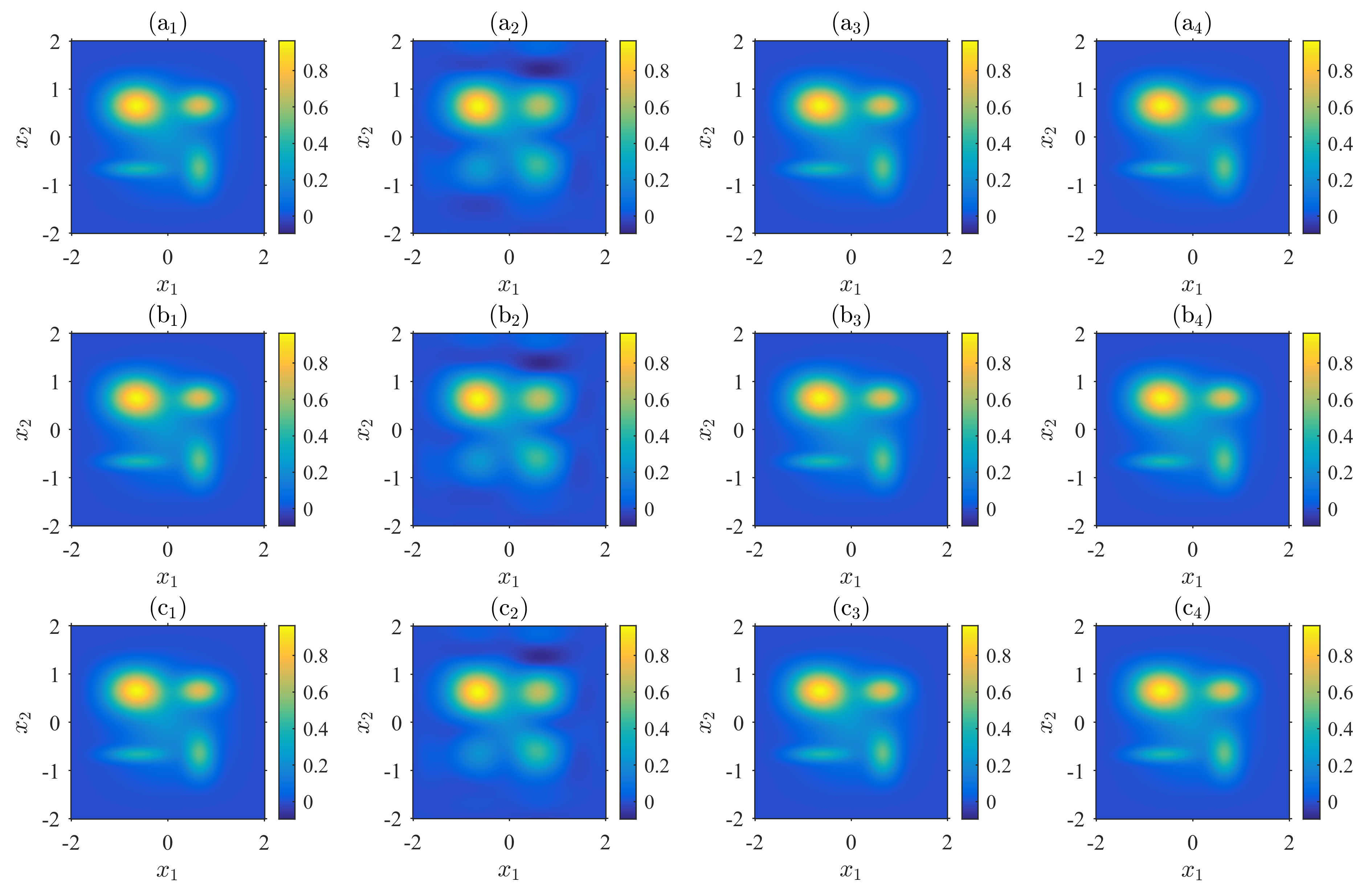}
		\caption{Reconstruction of the function $f$ for
			$\theta=\pi/4$. Rows correspond to $\mu=1,2,3$, respectively.
			The first column repeats the original function; the other columns show
			the reconstructions at $N=1,5,10$.}
		\label{fig:sobolev-reconstructions}
	\end{figure}

	The purpose of this section is to verify the effectiveness of the
	approximation formula in \eqref{sobolevseries}. Without loss of generality, we take $\theta=\pi/4$ and $\mu\in\{1, 2, 3\}$.   As an example, we consider a function $f$, which is defined as 
	$$
	f(x_1,x_2)=\sum_{j=1}^{5}c_j
	B_4\!\left(\frac{x_1-a_j}{r_j}\right)
	B_4\!\left(\frac{x_2-b_j}{s_j}\right),
	$$
	where 
	$$
	B_4(t)=
	\begin{cases}
		\dfrac{2}{3}-t^2+\dfrac{|t|^3}{2}, & |t|<1,\\[1mm]
		\dfrac{(2-|t|)^3}{6}, & 1\le |t|<2,\\[1mm]
		0, & |t|\ge2,
	\end{cases}
	$$ and the parameters $c_j,a_j,b_j,r_j,s_j$ are given in Table \ref{tab:sobolev-phantom}. A pseudocolor plot of $f$, whose support is
	$\operatorname{supp}f=[-2,2]^2$, is shown in the first column of
	Fig.\ref{fig:sobolev-reconstructions}.
	By \cite[Example 1.3]{HanShen}, $\nu_2(B_{4})=7/2$.
	Therefore, $\nu_2(f)=\nu_2(B_{4})=7/2>2$. Since $f$
	is compactly supported, we use the 16-point Gauss–Legendre quadrature rule to compute   the values 
	$\{V_{\theta,\mu}(2^{-N}\mathbf{k})\}$. The reconstruction version $f_{N}^{\#}$ is defined in  \eqref{sobolevseries},
	and the relative reconstruction error is defined as 
	\[
	\hbox{RES}:=\left(
	\frac{
		\sum_{\bm{x}_{i}\in \mathcal{G}}
		|f(\bm{x}_{i})-f_{N}^{\#}(\bm{x}_{i})|^2
	}{
		\sum_{\bm{x}_{i}\in \mathcal{G}}|f(\bm{x}_{i})|^2
	}
	\right)^{1/2}
	\]
	where 
	$\mathcal{G}=\{(-2+4m/10^3, -2+4n/10^3): m, n=0, \ldots, 10^3-1\}$.
	Table \ref{tab:sobolev-errors} reports the RES values, while
	Fig.\ref{fig:sobolev-reconstructions} displays the reconstructions
	$f^{\#}_{N}$ for different choices of $N$ and $\mu$.
	Both suggest that $f^{\#}_{N}$ converges to $f$ as $N\to\infty$.

	\begin{rem}
		There are three sources  of errors in the these numerical simulations: (1) the truncation error of the series in $\eqref{sobolevseries}$, where the summation is restricted to $\|2^{-N}\mathbf{k}\|_\infty\le L$, with $L=26,14,10$ corresponding to $\mu=1,2,3$, respectively; (2) Gauss–Legendre quadrature errors in computing the AV-LT samples $[V_{\theta,\mu}f](2^{-N}\mathbf{k})$; (3)  the approximate computation of the reconstruction kernel $\Phi$ using the FFT, as in subsection \ref{sec:numerical} \eqref{FFTjisuanAVL}. In the presence of these three sources of error, $f^{\#}_{N}$ stably approximates $f$ as the scale level $N$ increases.
	\end{rem}

	\begin{table}[!htbp]
		\centering
		\caption{Parameters of the test function.}
		\label{tab:sobolev-phantom}
		\renewcommand{\arraystretch}{1.20}
		\begin{tabular*}{0.85\textwidth}{@{\extracolsep{\fill}}cccccc@{}}
			\toprule
			$j$ & $c_j$ & $a_j$ & $b_j$ & $r_j$ & $s_j$ \\
			\midrule
			$1$ & $1/2$ & $0$ & $0$ & $1$ & $1$ \\
			$2$ & $2$ & $-2/3$ & $2/3$ & $2/3$ & $1/2$ \\
			$3$ & $3/2$ & $2/3$ & $2/3$ & $1/2$ & $1/3$ \\
			$4$ & $1$ & $2/3$ & $-2/3$ & $1/3$ & $1/2$ \\
			$5$ & $3/4$ & $-2/3$ & $-2/3$ & $2/3$ & $1/6$ \\
			\bottomrule
		\end{tabular*}
	\end{table}

	\begin{table}[!htbp]
		\centering
		\caption{Relative reconstruction error vs.  $(N, \mu)$.}
		\label{tab:sobolev-errors}
		\renewcommand{\arraystretch}{1.10}
		\begin{tabular*}{0.85\textwidth}{@{\extracolsep{\fill}}cccc@{}}
			\specialrule{\heavyrulewidth}{0pt}{0pt}
			\diagbox[width=4.6em,height=2.3em,linewidth=0.45pt,trim=l,innerleftsep=12pt,innerrightsep=12pt]{$N$}{$\mu$}
			& $1$ & $2$ & $3$ \\
			\specialrule{\lightrulewidth}{0pt}{\belowrulesep}
			$1$ & 1.5202e-1 & 1.5512e-1 & 1.5940e-1 \\
			$5$ & 1.8108e-5 & 1.8117e-5 & 1.8126e-5 \\
			$10$ & 2.0549e-10 & 1.9818e-10 & 1.9621e-10 \\
			\bottomrule
		\end{tabular*}
	\end{table}

	\FloatBarrier
	
	\appendix
	\section{Technical proofs for the measure \texorpdfstring{$\sigma_s$}{}}
	\label{Proof-measure}
	\subsection{Proof of Proposition \ref{def:measure} } \label{proofofpropositionmearue}
	
	Let $\mu>0$, and let $\theta\in(0,\pi/2)$ be fixed. For each $s\in\{+,-\}$ define
	\[
	\sigma_s(E):=\int_{0}^{\infty}\chi_E\bigl(-r\mathbf d_s(\theta)\bigr)e^{-\mu r}dr,
	\]
	where $E\subseteq\mathbb R^2$ is a Borel set.
	
	First we prove that $\sigma_s(E)$ is well-defined. Since the map
	$r\mapsto -r\mathbf d_s(\theta)$ is continuous from $[0,\infty)$ to $\mathbb R^2$, the function
	$r\mapsto \chi_E(-r\mathbf d_s(\theta))$ is Borel measurable on $[0,\infty)$ for every Borel set $E\subset\mathbb R^2$.
	Hence the integrand is measurable and nonnegative, so the integral is well-defined and lies in $[0,\infty]$.
	
	Clearly $\sigma_s(\varnothing)=0$. Let $\{E_n\}_{n\ge1}$ be pairwise disjoint Borel sets.
	Then for every $x\in\mathbb R^2$, $\chi_{\cup_{n=1}^\infty E_n}(x)=\sum_{n=1}^\infty \chi_{E_n}(x)$.
	Hence for every $r\ge0$, $\chi_{\cup_{n=1}^\infty E_n}\bigl(-r\mathbf d_s(\theta)\bigr)=\sum_{n=1}^\infty \chi_{E_n}\bigl(-r\mathbf d_s(\theta)\bigr)$.
	
	Define $g_N(r):=\sum_{n=1}^N \chi_{E_n}(-r\mathbf d_s(\theta))$.
	Then for every $r>0$ we have $g_N(r)\le g_{N+1}(r)$ and $\lim_{N\to\infty} g_N(r) = g(r)$, where
	$g(r)=\chi_{\cup_{n=1}^\infty E_n}(-r\mathbf d_s(\theta))$.
	By the monotone convergence theorem,
	$$
	\begin{aligned}
		\sigma_s\Bigl(\bigcup_{n=1}^\infty E_n\Bigr)
		&=\int_0^\infty \chi_{\cup_{n=1}^\infty E_n}(-r\mathbf d_s(\theta))e^{-\mu r}dr
		=\int_0^\infty \Bigl(\sum_{n=1}^\infty \chi_{E_n}(-r\mathbf d_s(\theta))\Bigr)e^{-\mu r}dr\\
		&=\lim_{N\to\infty}\int_0^\infty g_N(r)e^{-\mu r}dr
		=\lim_{N\to\infty}\sum_{n=1}^N\int_0^\infty \chi_{E_n}(-r\mathbf d_s(\theta))e^{-\mu r}dr\\
		&=\sum_{n=1}^\infty \sigma_s(E_n).
	\end{aligned}
	$$
	Thus $\sigma_s$ is a Borel measure on $\mathbb R^2$.
	
	Finally, a sum of two Borel measures is again a Borel measure.  This proves the proposition.

	\section{Technical proofs for the AV-LT characterizations}\label{app:characterization-proofs}
	This appendix collects the auxiliary estimates and longer constructions used in Section \ref{kehua}.
	
	\subsection{Proof of Lemma \ref{lem:modulation-estimate}}\label{Sobolev-nonsmoothing}
	It follows from $\widehat{\psi_N}(\boldsymbol{\xi})=\widehat{\psi}(\boldsymbol{\xi}-N\mathbf e)$
	that 
	\begin{equation}\label{eq:norm-identity}
		\|\psi_N\|_{H^\lambda(\mathbb R^2)}^2
		=\frac1{(2\pi)^2}
		\int_{\mathbb R^2}
		(1+\|\mathbf{\xi}+N\mathbf{e}\|_2^2)^\lambda
		|\widehat\psi(\mathbf{\xi})|^2d\mathbf{\xi}.
	\end{equation}
	
	We first prove the upper bound. By Peetre's inequality \cite[Proposition 3.3.31]{RuzhanskyTurunen2010},  for all $\mathbf{a},\mathbf{b}\in\mathbb R^2$ we have 
	\begin{equation}\label{eq:elementary}
		1+\|\mathbf{a}+\mathbf{b}\|_2^2
		\le
		2(1+\|\mathbf{a}\|_2^2)(1+\|\mathbf{b}\|_2^2).
	\end{equation}		
	Applying \eqref{eq:elementary} with $\mathbf{a}=\mathbf{\xi}$ and
	$\mathbf{b}=N\mathbf{e}$, we obtain
	\begin{equation}\label{eq:upper}
		(1+\|\mathbf{\xi}+N\mathbf{e}\|_2^2)^\lambda
		\le
		2^\lambda(1+N^2)^\lambda(1+\|\mathbf{\xi}\|_2^2)^\lambda.
	\end{equation}
	Combining \eqref{eq:norm-identity} and
	\eqref{eq:upper} gives
	\begin{align*}
		\|\psi_N\|_{H^\lambda(\mathbb R^2)}^2
		\lesssim
		(1+N^2)^\lambda
		\|\psi\|_{H^\lambda(\mathbb R^2)}^2
		\lesssim
		(1+N^2)^\lambda.    
	\end{align*}
	
	It remains to prove the lower bound. Applying  \eqref{eq:elementary} with
	$\mathbf{a}=\mathbf{\xi}+N\mathbf{e}$ and $\mathbf{b}=-\mathbf{\xi}$  yields
	\[
	1+N^2\le2(1+\|\mathbf{\xi}+N\mathbf{e}\|_2^2)(1+\|\mathbf{\xi}\|_2^2).
	\]
	Therefore, 
	\begin{equation}\label{eq:lower-weight-estimate}
		(1+\|\boldsymbol{\xi}+N\mathbf e\|_2^2)^\lambda
		\gtrsim
		(1+N^2)^\lambda
		(1+\|\boldsymbol{\xi}\|_2^2)^{-\lambda}.
	\end{equation}
	Substituting \eqref{eq:lower-weight-estimate} into
	\eqref{eq:norm-identity}, we get
	\[
	\|\psi_N\|_{H^\lambda(\mathbb R^2)}^2\gtrsim(1+N^2)^\lambda\int_{\mathbb R^2}(1+\|\boldsymbol{\xi}\|_2^2)^{-\lambda}|\widehat\psi(\boldsymbol{\xi})|^2d\boldsymbol{\xi},
	\]
	where the  integral is finite since
	$(1+\|\mathbf{\xi}\|_2^2)^{-\lambda}\le 1$ and
	$\widehat\psi\in L^2(\mathbb R^2)$.
	It is also strictly positive because $\psi\not\equiv0$ and
	$(1+\|\mathbf{\xi}\|_2^2)^{-\lambda}$ is positive on $\mathbb R^2$. Thus $\|\psi_N\|_{H^\lambda(\mathbb R^2)}\gtrsim(1+N^2)^{\lambda/2}$.
	The proof is completed.

	\subsection{Proof of Example \ref{ex:strict-sobolev} }\label{app:sobolev-details}
	The proof proceeds in two steps.

	\smallskip
	\underline{\textbf{Step 1}: Prove that $\nu_2(V_{\theta,\mu}f)\ge 5/2$}.
	By \eqref{Sobolevguanghuadingyi}, it suffices to prove that
	$V_{\theta,\mu}f\in H^s(\mathbb R^2)$ for every $0\le s<5/2$.
	Since $H^1(\mathbb R^2)\subseteq H^s(\mathbb R^2)$ for $0\le s<1$,
	this reduces to proving that $V_{\theta,\mu}f\in H^s(\mathbb R^2)$
	for every $1\le s<5/2$.
	By \eqref{eq:Fourier identity}, we therefore need to establish that
	\begin{equation}\label{eq:example-sobolev-integral}
		\|V_{\theta,\mu}f\|_{H^s(\mathbb R^2)}^2
		=\frac{1}{(2\pi)^2}\int_{\mathbb R^2}
		(1+\|\boldsymbol\xi\|_2^2)^s
		|m_{\theta,\mu}(\boldsymbol\xi)|^2
		|\widehat f(\boldsymbol\xi)|^2d\boldsymbol\xi
		<\infty
	\end{equation}
	for every $1\le s<5/2$.
	To this end, we split the integral over $C_\delta$ and its
	complement $C_\delta^c$, where
	\[
	C_\delta
	:=
	\left\{\boldsymbol\xi=(\xi_1,\xi_2)\in\mathbb R^2:
	|\xi_2|\le\delta|\xi_1|\right\}
	\cup
	\left\{\boldsymbol\xi=(\xi_1,\xi_2)\in\mathbb R^2:
	|\xi_1|\le\delta|\xi_2|\right\},
	\]
	with $0<\delta<\frac12\min\{\tan\theta,\cot\theta\}$.
	Then \eqref{eq:example-sobolev-integral} can be rewritten as
	$$
	\begin{aligned}
		\|V_{\theta,\mu}f\|_{H^s(\mathbb R^2)}^2
		&=\frac{1}{(2\pi)^2}\Biggl(
		\underbrace{\int_{C_\delta}
			(1+\|\boldsymbol\xi\|_2^2)^s
			|m_{\theta,\mu}(\boldsymbol\xi)|^2
			|\widehat f(\boldsymbol\xi)|^2d\boldsymbol\xi}_{I_1}\\
		&\qquad\qquad+
		\underbrace{\int_{C_\delta^c}
			(1+\|\boldsymbol\xi\|_2^2)^s
			|m_{\theta,\mu}(\boldsymbol\xi)|^2
			|\widehat f(\boldsymbol\xi)|^2d\boldsymbol\xi}_{I_2}
		\Biggr).
	\end{aligned}
	$$
	It remains to prove that $I_1<\infty$ and $I_2<\infty$.
	
	\underline{We first  estimate  $I_1$}.  Recall that $0<\delta<\frac12\min\{\tan\theta,\cot\theta\}$. For every $\boldsymbol{\xi} \in C_\delta$,
	if $|\xi_2|\le\delta|\xi_1|$, then
	\begin{align}\label{diyigetiaojian}
		|\xi_2|\cos\theta
		\le\delta|\xi_1|\cos\theta
		\le\tfrac12|\xi_1|\sin\theta.
	\end{align}
	Similarly, if $|\xi_1|\le\delta|\xi_2|$, then
	\begin{align}\label{diergetiaojian}
		|\xi_1|\sin\theta
		\le\delta|\xi_2|\sin\theta
		\le\tfrac12|\xi_2|\cos\theta.
	\end{align}
	Thus, 
	\begin{equation}\label{example-dkesei}
		\begin{aligned}
			|\mathbf d_\pm(\theta)\cdot\boldsymbol\xi|
			&\ge \bigl||\xi_1|\sin\theta-|\xi_2|\cos\theta\bigr| \\
			&\ge \frac12
			\max\bigl\{|\xi_1|\sin\theta,|\xi_2|\cos\theta\bigr\}
			\quad  (\hbox{from} \ \  \eqref{diyigetiaojian} \  \hbox{and} \  \eqref{diergetiaojian}) 
			\\
			&\ge \frac{\min\{\sin\theta,\cos\theta\}}{2}
			\max\{|\xi_1|,|\xi_2|\} \\
			&\ge \frac{\min\{\sin\theta,\cos\theta\}}{2\sqrt2}
			\|\boldsymbol\xi\|_2.
		\end{aligned}
	\end{equation}
	So, by  \eqref{eq:m-theta-mu} and \eqref{example-dkesei} we have 
	$$
	|m_{\theta,\mu}(\boldsymbol\xi)|
	\le\frac{1}{\sqrt{\mu^2+|\mathbf d_+(\theta)\cdot\boldsymbol\xi|^2}}
	+\frac{1}{\sqrt{\mu^2+|\mathbf d_-(\theta)\cdot\boldsymbol\xi|^2}}
	\lesssim(1+\|\boldsymbol\xi\|_2^2)^{-1/2},
	\qquad\boldsymbol\xi\in C_\delta.
	$$
	Since $0\le s-1<3/2=\nu_2(f)$, we obtain
	\begin{align}\label{gujiI1}
		I_1\lesssim\int_{\mathbb R^2}
		(1+\|\boldsymbol\xi\|_2^2)^{s-1}
		|\widehat f(\boldsymbol\xi)|^2d\boldsymbol\xi
		=(2\pi)^2\|f\|_{H^{s-1}(\mathbb R^2)}^2<\infty.
	\end{align}
	
	\underline{Next we estimate  $I_2$}. 
	By the convolution theorem, we have 
	\begin{equation}\label{example-f-Fourier}
		\widehat f(\xi_1,\xi_2)
		=[\widehat{\chi_{(0,1]}}(\xi_1)]^2
		[\widehat{\chi_{(0,1]}}(\xi_2)]^2
		=e^{-\mathrm i(\xi_1+\xi_2)}
		\prod_{j=1}^2\left(\frac{\sin(\xi_j/2)}{\xi_j/2}\right)^2.
	\end{equation}
	The elementary estimate
	$$
	\left|\frac{\sin(t/2)}{t/2}\right|
	\le\min\left\{1,\frac2{|t|}\right\}
	\lesssim(1+|t|)^{-1},
	\qquad t\ne0,
	$$
	yields $|\widehat f(\xi_1,\xi_2)|\lesssim(1+|\xi_1|)^{-2}(1+|\xi_2|)^{-2}$.
	Moreover, we have $|\xi_1|\asymp|\xi_2|\asymp\|\boldsymbol\xi\|_2$ for every $\boldsymbol\xi \in C_\delta^c$. Therefore, we arrive at $|\widehat f(\boldsymbol\xi)|^2\lesssim(1+\|\boldsymbol\xi\|_2^2)^{-4}$.
	
	On the other hand, using $|m_{\theta,\mu}(\boldsymbol\xi)|\le 2/\mu$, we obtain
	\begin{align}\label{gujiI2}
		I_2\le
		\frac{4}{\mu^2}
		\int_{C_\delta^c}
		(1+\|\boldsymbol\xi\|_2^2)^s
		|\widehat f(\boldsymbol\xi)|^2d\boldsymbol\xi
		\lesssim
		\int_{\mathbb R^2}
		(1+\|\boldsymbol\xi\|_2^2)^{s-4}d\boldsymbol\xi,
	\end{align}
	where the last integral converges for $s<5/2<3$. 
	Combining \eqref{gujiI1} and \eqref{gujiI2}, we have 
	$V_{\theta,\mu}f\in H^s(\mathbb R^2)$ for every $1\le s<5/2$. Hence, $\nu_2(V_{\theta,\mu}f)\ge5/2$.

	\smallskip
	\underline{\textbf{Step 2}: Prove   $\nu_2(V_{\theta,\mu}f)\le 5/2$}.
	To obtain $\nu_2(V_{\theta,\mu}f)\le5/2$, it suffices to prove that
	$V_{\theta,\mu}f\notin H^s(\mathbb R^2)$ for every $s>5/2$.
	Suppose, to the contrary, that $V_{\theta,\mu}f\in H^s(\mathbb R^2)$
	for some $s>5/2$. We will show that this contradicts $\nu_2(f)=3/2$.
	
	Since $\bigl|\mu-\mathrm i\mathbf d_\pm(\theta)\cdot\boldsymbol\xi\bigr|^2
	=\mu^2+|\mathbf d_\pm(\theta)\cdot\boldsymbol\xi|^2
	\le\mu^2+\|\boldsymbol\xi\|_2^2$, \eqref{eq:m-theta-mu} yields that
	\begin{equation}\label{example-m-lowbound}
		|m_{\theta,\mu}(\xi_1,\xi_2)|^2
		=\frac{4(\mu^2+\cos^2\theta\xi_2^2)}
		{\bigl|\mu-\mathrm i\mathbf d_+(\theta)\cdot\boldsymbol\xi\bigr|^2
			\bigl|\mu-\mathrm i\mathbf d_-(\theta)\cdot\boldsymbol\xi\bigr|^2}
		\gtrsim\frac{1+\xi_2^2}
		{(1+\|\boldsymbol\xi\|_2^2)^2}.
	\end{equation}
	Similarly, interchanging $\xi_1$ and $\xi_2$ and adding \eqref{example-m-lowbound}, we obtain
	$$
	|m_{\theta,\mu}(\xi_1,\xi_2)|^2
	+|m_{\theta,\mu}(\xi_2,\xi_1)|^2
	\gtrsim\frac{2+\xi_1^2+\xi_2^2}
	{(1+\|\boldsymbol\xi\|_2^2)^2}
	\gtrsim\frac1{1+\|\boldsymbol\xi\|_2^2}.
	$$ 
	Moreover, from \eqref{example-f-Fourier}, we have $\widehat f(\xi_1,\xi_2)=\widehat f(\xi_2,\xi_1)$. By symmetry, we obtain
	$$
	\begin{aligned}
		\|V_{\theta,\mu}f\|_{H^s(\mathbb R^2)}^2
		&=\frac1{2(2\pi)^2}\int_{\mathbb R^2}
		(1+\|\boldsymbol\xi\|_2^2)^s
		\Bigl(|m_{\theta,\mu}(\xi_1,\xi_2)|^2
		+|m_{\theta,\mu}(\xi_2,\xi_1)|^2\Bigr)
		|\widehat f(\boldsymbol\xi)|^2d\boldsymbol\xi\\
		&\gtrsim\int_{\mathbb R^2}
		(1+\|\boldsymbol\xi\|_2^2)^{s-1}
		|\widehat f(\boldsymbol\xi)|^2d\boldsymbol\xi
		=(2\pi)^2\|f\|_{H^{s-1}(\mathbb R^2)}^2.
	\end{aligned}
	$$
	
	By the hypothesis, if $V_{\theta,\mu}f\in H^s(\mathbb R^2)$
	for some $s>5/2$, it follows that $f\in H^{s-1}(\mathbb R^2)$. Since $s-1>3/2$, this contradicts $\nu_2(f)=3/2$. Hence $V_{\theta,\mu}f\notin H^s(\mathbb R^2)$ for every $s>5/2$, and $\nu_2(V_{\theta,\mu}f)\le5/2$.
	Together with step 1, this proves $\nu_2(V_{\theta,\mu}f)=5/2$.

	\subsection{Proof of Proposition \ref{compact-example}}\label{app:compact-example-details}
	Clearly,  $f$ is compactly supported.
	From Theorem \ref{L2-boundedness}, $V_{\theta,\mu}f\in L^{2}(\mathbb{R}^{2})$.
	We next prove that it is also  compactly supported. 
	For $\boldsymbol{\xi},\mathbf{a}\in\mathbb{R}^2$, we have 
	$\widehat{g(\cdot+\mathbf{a})}(\boldsymbol{\xi})
	=
	e^{\mathrm{i}\boldsymbol{\xi}\cdot\mathbf{a}}\widehat{g}(\boldsymbol{\xi}),$
	from which  we arrive at 
	\[
	\begin{aligned}
		&\widehat{f}(\boldsymbol{\xi})=\\
		&-\frac{1}{4}\biggl(
		e^{-2\mu}e^{\mathrm{i}\boldsymbol{\xi}\cdot(\mathbf{d}_+(\theta)+\mathbf{d}_-(\theta))}
		- e^{\mathrm{i}\boldsymbol{\xi}\cdot(\mathbf{d}_+(\theta)-\mathbf{d}_-(\theta))}
		- e^{\mathrm{i}\boldsymbol{\xi}\cdot(-\mathbf{d}_+(\theta)+\mathbf{d}_-(\theta))}
		+ e^{2\mu}e^{\mathrm{i}\boldsymbol{\xi}\cdot(-\mathbf{d}_+(\theta)-\mathbf{d}_-(\theta))}
		\biggr)\widehat{g}(\boldsymbol{\xi}).
	\end{aligned}
	\]
	On the other hand, using 
	$\sin z=
	\frac{e^{\mathrm{i}z}-e^{-\mathrm{i}z}}{2\mathrm{i}}$
	with 
	$z=\boldsymbol{\xi}\cdot\mathbf{d}_{+}(\theta)+\mathrm{i}\mu$ and $
	z=\boldsymbol{\xi}\cdot\mathbf{d}_{-}(\theta)+\mathrm{i}\mu,$
	we obtain
	\begin{align}\label{FTEX}
		\widehat{f}(\boldsymbol{\xi})
		=
		\sin\bigl(\boldsymbol{\xi}\cdot\mathbf{d}_{+}(\theta)+\mathrm{i}\mu\bigr)
		\sin\bigl(\boldsymbol{\xi}\cdot\mathbf{d}_{-}(\theta)+\mathrm{i}\mu\bigr)
		\widehat{g}(\boldsymbol{\xi}),
		\ 
		\boldsymbol{\xi}\in\mathbb{R}^2. 
	\end{align}

	
	Recall that  $f$ is compactly supported. Then $\widehat f$  can extend to an entire function on $\mathbb C^2$. Hence \eqref{FTEX} extends to all $\mathbf\zeta\in\mathbb C^2$:
	\[
	\widehat f(\mathbf\zeta)
	=\sin\bigl(\mathbf\zeta\cdot \mathbf d_{+}(\theta)+\mathrm i\mu\bigr)
	\sin\bigl(\mathbf\zeta\cdot \mathbf d_{-}(\theta)+\mathrm i\mu\bigr)
	\widehat g(\mathbf\zeta).
	\]
	Now it is sufficient to verify Theorem \ref{thm:compact-support-V} (2). We fix $s\in\{+,-\}$ and $\omega\in\mathbb R$, and set
	$\mathbf\zeta=\omega\mathbf d_s^{\perp}(\theta)-\mathrm i\mu\mathbf d_s(\theta)\in\mathbb C^2$. Then $\mathbf\zeta\cdot \mathbf d_s(\theta)=-\mathrm i\mu$, and therefore $\sin\bigl(\mathbf\zeta\cdot \mathbf d_s(\theta)+\mathrm i\mu\bigr)=0$, which implies $\widehat f(\mathbf\zeta)=0$. This  completes the proof. 
	
	\section{Technical proofs for the reconstruction theory}\label{app:technical-proofs}

	\subsection{Proof of Lemma \ref{lem:min-m}}\label{app:riesz-details}
	By \eqref{eq:m-theta-mu}, we have
	\[
	|m_{\theta,\mu}(\xi_1,\xi_2)|^2
	=
	\frac{4(\mu^2+\xi_2^2\cos^2\theta)}
	{\bigl(\mu^2+(\xi_1\sin\theta+\xi_2\cos\theta)^2\bigr)
		\bigl(\mu^2+(\xi_2\cos\theta-\xi_1\sin\theta)^2\bigr)}.
	\]
	Inspired  by 
	\[
	\bigl|m_{\theta,\mu}(\xi_1,\xi_2)\bigr|^2
	=
	\bigl|m_{\theta,\mu}(-\xi_1,\xi_2)\bigr|^2
	=
	\bigl|m_{\theta,\mu}(\xi_1,-\xi_2)\bigr|^2,
	\]
	we define  $X:=|\xi_1|$ and  $Y:=|\xi_2|$. By $\eqref{def:ho,R}$, we have $X\in[\rho_1,R_1]$ and  $Y\in[\rho_2,R_2]$.
	Setting $u:=X^2\sin^2\theta$ and  $v:=Y^2\cos^2\theta$,   we just need to consider the maximum and minimum of 
	\[
	G(u,v)
	:=
	\frac{4(\mu^2+v)}
	{(\mu^2+u-v)^2+4\mu^2v}
	\]
	on the rectangle
	$[U_0,U_1]\times[V_0,V_1]$.
	
	For each fixed $v\in[V_0,V_1]$, direct differentiation yields
	\begin{equation}\label{G-derivation}
		\frac{\partial G}{\partial u}(u,v)
		=
		-\frac{8(\mu^2+v)(\mu^2+u-v)}
		{\bigl((\mu^2+u-v)^2+4\mu^2v\bigr)^2}.
	\end{equation}
	
	Since $\mu>0$ and $v\geq 0$, the denominator of $\eqref{G-derivation}$ is positive and
	$\mu^2+v>0$. Hence the sign of $\partial G/\partial u$ is determined
	only by $-(\mu^2+u-v)$. Consequently, for each
	fixed $v\in[V_0,V_1]$, the minimum of $G(\cdot,v)$ over
	$[U_0,U_1]$ is attained at one of the endpoints, $u=U_0$ or $u=U_1$.
	It remains to minimize $G$ on the two sides $u=U_0$ and
	$u=U_1$. For fixed $u$, we have
	\[
	\frac{\partial G}{\partial v}(u,v)
	=
	\frac{4\bigl(u^2-v^2-\mu^4+4\mu^2u-2\mu^2v\bigr)}
	{\bigl((\mu^2+u-v)^2+4\mu^2v\bigr)^2}.
	\]
	Set $P_u(v):=u^2-v^2-\mu^4+4\mu^2u-2\mu^2v$. Then $P_u'(v)=-2v-2\mu^2<0$, and hence $P_u$ is strictly decreasing on $[V_0,V_1]$. It follows that, on
	each side $u=U_0$ and $u=U_1$, the minimum is attained at one of the
	endpoints, $v=V_0$ or $v=V_1$. Therefore, $\min_{\boldsymbol{\xi}\in \Omega}|m_{\theta,\mu}(\boldsymbol{\xi})|^2$ is attained at one of the four vertices
	\[
	(U_0,V_0),\quad (U_0,V_1),\quad (U_1,V_0),\quad (U_1,V_1).
	\]
	
	We now consider $\max_{\boldsymbol{\xi}\in \Omega}|m_{\theta,\mu}(\boldsymbol{\xi})|^2$. Suppose that $G$ has an interior critical point $(u_0,v_0)$. Then $\frac{\partial G}{\partial u}(u_0,v_0)=0$ and $\frac{\partial G}{\partial v}(u_0,v_0)=0$. The first equality implies $u_0=v_0-\mu^2$. The second equality would force the numerator of $\partial G/\partial v$ at $(u_0,v_0)$ to vanish. However, substituting $u_0=v_0-\mu^2$ gives
	\[
	u_0^2-v_0^2-\mu^4+4\mu^2u_0-2\mu^2v_0
	=-4\mu^4\neq 0,
	\]
	which contradicts $\frac{\partial G}{\partial v}(u_0,v_0)=0$. Hence $G$ has no interior critical point.
	Therefore, $\max_{\boldsymbol{\xi}\in \Omega}|m_{\theta,\mu}(\boldsymbol{\xi})|^2$ is attained on the boundary of $[U_0,U_1]\times[V_0,V_1]$.
	The proof is completed.
	\subsection{Proof of Theorem \ref{th:Riesz-bounds}}\label{proof-riesz-basis}
	For any finite sequence $c=\{c_{\boldsymbol{k}}\}_{{\boldsymbol{k}}\in\mathbb Z^2}$, define
	\[
	f_c(\mathbf{x})
	:=
	\sum_{k\in\mathbb Z^2}
	c_k\Phi(\mathbf{x}-\lambda_{\boldsymbol{k}}),
	\quad \mathbf{x}\in\mathbb R^2.
	\]
	Then $\widehat{f_c}(\mathbf{\xi})=\widehat{\Phi}(\mathbf{\xi})F_c(\mathbf{\xi})$, where $F_c(\mathbf{\xi}):=\sum_{{\boldsymbol{k}}\in\mathbb Z^2}c_{\boldsymbol{k}} e^{-\mathrm{i}\lambda_{\boldsymbol{k}}\cdot\mathbf{\xi}}$. From this and \eqref{Fourierexpression}, $f_c\in PW_\Omega$. By Parseval's
	identity,
	\[\|F_c\|_{L^2(\Omega)}^2=\sum_{\boldsymbol{k}\in\mathbb Z^2}|c_{\boldsymbol{k}}|^2,
	\]
	where $L^2(\Omega)$ is equipped with the normalized inner product in \eqref{neijikongjian}.
	
	On the other hand, Plancherel's theorem implies that
	\begin{align}\label{jifenbiaoshi}
		\|f_c\|_{L^2(\mathbb R^2)}^2=\frac{1}{(2\pi)^2}\int_{\mathbb R^2}
		|\widehat{f_c}(\boldsymbol{\xi})|^2d\boldsymbol{\xi} =
		\frac{1}{(2\pi)^2}\int_\Omega|\widehat{\Phi}(\boldsymbol{\xi})|^2
		|F_c(\boldsymbol{\xi})|^2d\boldsymbol{\xi}.
	\end{align}
	Define
	\begin{equation}\label{def:H}
		H(\boldsymbol{\xi}):=|\widehat{\Phi}(\boldsymbol{\xi})|^2=h_1^2h_2^2
		\frac{\chi_\Omega(\boldsymbol{\xi})}
		{|m_{\theta,\mu}(\boldsymbol{\xi})|^2},
		\quad \boldsymbol{\xi}\in\mathbb R^2.
	\end{equation}
	By Remark \ref{rem:m-properties}, the function $H$ attains its minimum and maximum on $\Omega$. Define 
	\[
	A_{\mathrm{opt}}
	:=
	\frac{1}{h_1h_2}
	\min_{\boldsymbol{\xi}\in\Omega}H(\boldsymbol{\xi}),
	\quad
	B_{\mathrm{opt}}
	:=
	\frac{1}{h_1h_2}
	\max_{\boldsymbol{\xi}\in\Omega}H(\boldsymbol{\xi}).
	\]
	Now it follows from \eqref{jifenbiaoshi} that 
	\begin{align}\label{you}
		\|f_c\|_{L^2(\mathbb R^2)}^2\le
		\left(\max_{\boldsymbol{\xi}\in\Omega}H(\boldsymbol{\xi})\right)
		\frac{1}{(2\pi)^2}
		\int_\Omega |F_c(\boldsymbol{\xi})|^2d\boldsymbol{\xi}=
		\left(\max_{\boldsymbol{\xi}\in\Omega}H(\boldsymbol{\xi})\right)
		\frac{1}{h_1h_2}
		\|F_c\|_{L^2(\Omega)}^2,
	\end{align}
	and 
	\begin{align}\label{zuo}
		\|f_c\|_{L^2(\mathbb{R}^2)}^2
		\ge
		\left(\min_{\boldsymbol{\xi}\in\Omega}H(\boldsymbol{\xi})\right)
		\frac{1}{(2\pi)^2}
		\int_\Omega |F_c(\boldsymbol{\xi})|^2d\boldsymbol{\xi}
		=\left(\min_{\boldsymbol{\xi}\in\Omega}H(\boldsymbol{\xi})\right)
		\frac{1}{h_1h_2}
		\|F_c\|_{L^2(\Omega)}^2.
	\end{align}

	The above bounds are sharp, namely, 
	the  optimal Riesz lower
	and upper bounds are $A_{\mathrm{opt}}$
	and $B_{\mathrm{opt}}$, respectively, such that  
	\[
	A_{\mathrm{opt}}\|F_c\|_{L^2(\Omega)}^2
	\le
	\|f_c\|_{L^2(\mathbb R^2)}^2
	\le
	B_{\mathrm{opt}}\|F_c\|_{L^2(\Omega)}^2.
	\]
	By $\eqref{def:H}$, we get
	\[
	A_{\mathrm{opt}}
	=
	\frac{h_1h_2}
	{\displaystyle\max_{\xi\in\Omega}|m_{\theta,\mu}(\boldsymbol{\xi})|^2},
	\quad
	B_{\mathrm{opt}}
	=
	\frac{h_1h_2}
	{\displaystyle\min_{\xi\in\Omega}|m_{\theta,\mu}(\boldsymbol{\xi})|^2}.
	\]
	where  the maximum and minimum of $|m_{\theta,\mu}|^2$ on $\Omega$ is given  in Lemma \ref{lem:min-m}. This yields $A_{\mathrm{opt}}$ and $B_{\mathrm{opt}}$. 

	To prove the  completeness, let $f\in PW_{\Omega}$ be orthogonal to
	$\Phi(\cdot-\lambda_{\boldsymbol{k}})$ for all $\mathbf{k}\in\mathbb Z^2$. Then, by Plancherel's theorem,
	\[
	0=\frac{1}{(2\pi)^2}\int_{\Omega}\widehat f(\boldsymbol{\xi})
	\overline{\widehat{\Phi}(\boldsymbol{\xi})}e^{\mathrm{i}\lambda_{\boldsymbol{k}}\cdot\boldsymbol{\xi}}d\boldsymbol{\xi},
	\qquad \forall \mathbf{k}\in\mathbb Z^2.
	\]
	Thus $\widehat f\overline{\widehat{\Phi}}$ is orthogonal to all exponentials $\{e^{-\mathrm{i}{\lambda_{\boldsymbol{k}}}\cdot\boldsymbol{\xi}}\}_{\mathbf{k}\in\mathbb Z^2}$, which form an orthogonal basis of $L^2(\Omega)$. Hence $\widehat f(\boldsymbol{\xi})\overline{\widehat{\Phi}(\boldsymbol{\xi})}=0$ a.e. on $\Omega$. Since $\widehat{\Phi}\neq 0$ on $\Omega$, we get $\widehat f=0$ a.e. on $\Omega$. It follows that $f\equiv 0$. Therefore
	$\left\{\Phi(\cdot-\lambda_{\boldsymbol{k}})\right\}_{\boldsymbol{k}\in\mathbb Z^2}$ is complete in $PW_{\Omega}$, and the proof is completed.
	
	\subsection{Proof of Lemma \ref{lem:Phi-decay}}\label{app:truncation-details}
	For fixed $\xi_2\in[a_2,b_2]$, define
	\[
	I(x_1,\xi_2)
	:=
	\int_{a_1}^{b_1}
	a(\xi_1,\xi_2)e^{\mathrm{i}x_1\xi_1}d\xi_1.
	\]
	Then
	\[
	\Phi(x_1,x_2)
	=
	\frac{h_1h_2}{(2\pi)^2}
	\int_{a_2}^{b_2}
	I(x_1,\xi_2)e^{\mathrm{i}x_2\xi_2}d\xi_2.
	\]
	\underline{\textbf{Estimate of $I$}.}
	For $|x_1|\le 1$, we first estimate that
	\begin{align}\label{firstestimate}
		|I(x_1,\xi_2)|
		\le
		\frac{2\pi}{h_1}M_0^\Omega
		\le
		\frac{4\pi(M_0^\Omega+M_1^\Omega)}{h_1(1+|x_1|)}.
	\end{align}
	
	On the other hand, for $|x_1|>1$, integration by parts in $\xi_1$ yields
	\begin{align}\label{zhongjianguocheng}
		I(x_1,\xi_2)
		=
		\frac{
			a(b_1,\xi_2)e^{\mathrm{i}b_1x_1}
			-
			a(a_1,\xi_2)e^{\mathrm{i}a_1x_1}
		}{\mathrm{i}x_1}
		-
		\frac{1}{\mathrm{i}x_1}
		\int_{a_1}^{b_1}
		\partial_{\xi_1}a(\xi_1,\xi_2)e^{\mathrm{i}x_1\xi_1}d\xi_1,  
	\end{align}
	which leads to
	\begin{align}\label{secondestimate}
		|I(x_1,\xi_2)|
		\le
		\frac{2h_1M_0^\Omega+2\pi M_1^\Omega}{h_1|x_1|}
		\le
		C_\Omega
		\frac{M_0^\Omega+M_1^\Omega}{1+|x_1|}. 
	\end{align}
	Combining \eqref{firstestimate} and 
	\eqref{secondestimate}, we obtain
	\begin{align}\label{eq:IN-bound}
		|I(x_1,\xi_2)|
		\le
		C_\Omega
		\frac{M_0^\Omega+M_1^\Omega}{1+|x_1|}.
	\end{align}
	\underline{\textbf{Estimate of $\partial_{\xi_2}I$}}.	
	Through the similar procedures in 
	\eqref{firstestimate}, \eqref{zhongjianguocheng} and \eqref{secondestimate},
	we have 
	\begin{align}\label{eq:dIN-bound}
		|\partial_{\xi_2}I(x_1,\xi_2)|
		\le
		C_\Omega
		\frac{M_1^\Omega+M_2^\Omega}{1+|x_1|}.
	\end{align}
	\underline{\textbf{Estimate of $\Phi$}}. For  $|x_2|\le 1$,  it follows from  \eqref{eq:IN-bound} that 
	\[
	|\Phi(x_1,x_2)|
	\le
	C_\Omega
	\frac{M_0^\Omega+M_1^\Omega+M_2^\Omega}{(1+|x_1|)(1+|x_2|)}.
	\]
	For $|x_2|>1$, through the similar analysis in \eqref{zhongjianguocheng} we have 
	
	\[
	|\Phi(x_1,x_2)|
	\le
	C_\Omega
	\frac{M_0^\Omega+M_1^\Omega+M_2^\Omega}{(1+|x_1|)(1+|x_2|)}.
	\]
	
	The proof is completed.

	\subsection{Proof of Lemma \ref{lem:Phi-lower-bound}}\label{proof-fuzhu-zuiyou}
	By \eqref{TTCX} we have 
	\[
	\Phi(t,0)
	=
	\frac{h_1h_2}{(2\pi)^2}
	\int_{a_2}^{b_2}\int_{a_1}^{b_1}
	a(\xi_1,\eta)e^{\mathrm{i}t\xi_1}d\xi_1d\eta .
	\]
	Since $a\in C^1(\overline\Omega)$, integration by parts with respect to
	$\xi_1$ gives, for $t\ne0$,
	\[
	\begin{aligned}
		\Phi(t,0)
		&=
		\frac{h_1h_2}{(2\pi)^2\mathrm{i}t}
		\left[
		e^{\mathrm{i}b_1t}\int_{a_2}^{b_2}a(b_1,\eta)d\eta
		-
		e^{\mathrm{i}a_1t}\int_{a_2}^{b_2}a(a_1,\eta)d\eta
		\right]                                                     \\
		&\quad
		-
		\frac{h_1h_2}{(2\pi)^2\mathrm{i}t}
		\int_{a_2}^{b_2}\int_{a_1}^{b_1}
		\partial_{\xi_1}a(\xi_1,\eta)e^{\mathrm{i}t\xi_1} d\xi_1d\eta .
	\end{aligned}
	\]
	The last integral is $o(1)$ as $|t|\to\infty$ by 
	Riemann-Lebesgue lemma. Thus
	\[
	\Phi(t,0)
	=
	\frac{h_1h_2}{(2\pi)^2\mathrm{i}t}
	\left[
	e^{\mathrm{i}b_1t}\int_{a_2}^{b_2}a(b_1,\eta)d\eta
	-
	e^{\mathrm{i}a_1t}\int_{a_2}^{b_2}a(a_1,\eta)d\eta
	\right]
	+
	o\left(\frac1{|t|}\right).
	\]
	Since $h_1=2\pi/(b_1-a_1)$, we have $(b_1-a_1)t_n=2\pi\left(\frac12-n\right)$,
	and therefore $e^{\mathrm{i}(b_1-a_1)t_n}=-1$.
	Hence we have $e^{\mathrm{i}b_1t_n}$
	$=-e^{\mathrm{i}a_1t_n}$ and
	\[
	\Phi(t_n,0)
	=
	-\frac{h_1h_2e^{\mathrm{i}a_1t_n}}{(2\pi)^2\mathrm{i}t_n}
	\left[
	\int_{a_2}^{b_2}a(b_1,\eta)d\eta
	+
	\int_{a_2}^{b_2}a(a_1,\eta)d\eta
	\right]
	+
	o\left(\frac1{|t_n|}\right).
	\]
	Note that the coefficient in brackets is nonzero. Indeed, by Remark \ref{rem:m-properties} we have 
	\[
	\operatorname{Re}a(\boldsymbol\xi)
	=
	\frac{\mu}{2}
	+
	\frac{\mu\xi_1^2\sin^2\theta}
	{2(\mu^2+\xi_2^2\cos^2\theta)}
	>0,
	\qquad \boldsymbol\xi\in\Omega.
	\]
	Therefore,
	\[
	\operatorname{Re}
	\left[
	\int_{a_2}^{b_2}a(b_1,\eta)d\eta
	+
	\int_{a_2}^{b_2}a(a_1,\eta)d\eta
	\right]
	>0.
	\]
	Together with the relation $|t_n|\asymp n$, this implies that there exist a  constant $N_0\in\mathbb N$ such that $|\Phi(t_n,0)|\gtrsim 1/n$ whenever $n\ge N_0$. The lemma is proved.
	
	\subsection{Proof of Theorem \ref{thm:truncation}}\label{proof-zuiyou}
	By Theorem \ref{thm:sampling-reconstruction}, for every $f\in PW_\Omega$,
	\[
	f(\mathbf x)
	=
	\sum_{\mathbf k\in\mathbb Z^2}
	[V_{\theta,\mu}f](\lambda_{\mathbf k})
	\Phi(\mathbf x-\lambda_{\mathbf k}),
	\qquad \mathbf x\in\mathbb R^2.
	\]
	Denote 
	\[
	E_Rf(\mathbf x):=
	\sum_{\substack{\mathbf k\in\mathbb Z^2\\
			\|\mathbf x-\lambda_{\mathbf k}\|_2>R}}
	[V_{\theta,\mu}f](\lambda_{\mathbf k})
	\Phi(\mathbf x-\lambda_{\mathbf k}).
	\]
	Then 
	\begin{align}\label{yuxiangguji}
		&|E_Rf(\mathbf x)|
		\le
		\left(
		\sum_{\mathbf k\in\mathbb Z^2}
		|[V_{\theta,\mu}f](\lambda_{\mathbf k})|^2
		\right)^{1/2}
		\left(
		\sum_{\substack{\mathbf k\in\mathbb Z^2\\
				\|\mathbf x-\lambda_{\mathbf k}\|_2>R}}
		|\Phi(\mathbf x-\lambda_{\mathbf k})|^2
		\right)^{1/2} \notag                                      \\
		&\lesssim
		\|f\|_{L^2(\mathbb R^2)}
		\left(
		\sum_{\substack{\mathbf k\in\mathbb Z^2\\
				\|\mathbf x-\lambda_{\mathbf k}\|_2>R}}
		\frac{1}{(1+|x_1-h_1k_1|)^2(1+|x_2-h_2k_2|)^2}
		\right)^{1/2} (\hbox{from} \ \eqref{1111} \ \hbox{and} \ \eqref{shengchengyuanguji})      
	\end{align}
	Next we establish the further estimate for \eqref{yuxiangguji}.
	We need 
	$\sup_{y\in\mathbb R}\sum_{k\in\mathbb Z}\frac1{(1+|y-hk|)^2}< \infty$ 
	and 
	\[
	\sup_{y\in\mathbb R}\sum_{\substack{k\in\mathbb Z\\ |y-hk|>r}}\frac1{(1+|y-hk|)^2}\lesssim r^{-1},
	\]
	which is derived  from 	 $\sum_{k\geqslant n\geqslant 0}{\frac{1}{\left( 1+k \right) ^2}\lesssim \int_n^{\infty}{\frac{dx}{\left( 1+x \right) ^2}}}$.     
	Now let $r=R/\sqrt2$. Since
	\[
	\{\mathbf k:\|\mathbf x-\lambda_{\mathbf k}\|_2>R\}
	\subset
	\{\mathbf k:|x_1-h_1k_1|>r\}
	\cup
	\{\mathbf k:|x_2-h_2k_2|>r\},
	\]
	we get
	\[
	\begin{aligned}
		&\sum_{\substack{\mathbf k\in\mathbb Z^2\\
				\|\mathbf x-\lambda_{\mathbf k}\|_2>R}}
		\frac1{(1+|x_1-h_1k_1|)^2(1+|x_2-h_2k_2|)^2}                 \\
		&\le
		\sum_{\substack{k_1\in\mathbb Z\\ |x_1-h_1k_1|>r}}
		\frac1{(1+|x_1-h_1k_1|)^2}
		\sum_{k_2\in\mathbb Z}
		\frac1{(1+|x_2-h_2k_2|)^2}                                               \\
		&\quad+
		\sum_{k_1\in\mathbb Z}
		\frac1{(1+|x_1-h_1k_1|)^2}
		\sum_{\substack{k_2\in\mathbb Z\\ |x_2-h_2k_2|>r}}
		\frac1{(1+|x_2-h_2k_2|)^2}                                                     \\
		&\lesssim R^{-1}.
	\end{aligned}
	\]
	We obtain
	$\sup_{\mathbf{x}\in\mathbb{R}^2}\bigl|f(\mathbf{x})-S_Rf(\mathbf{x})\bigr|\lesssim R^{-1/2}$.
	
	It remains to prove the convergence rate $1/2$ is optimal. We will construct a function $f\in PW_{\Omega}$ such that, for every $\delta>0$, $\sup_{\mathbf x\in\mathbb R^2}|f(\mathbf x)-S_Rf(\mathbf x)|\ne O(R^{-1/2-\delta})$ as $R\rightarrow \infty$.
	Let $\mathbf x^0=\left(\frac{h_1}{2},0\right), t_n=h_1\left(\frac12-n\right), n\in\mathbb N$. Then $\mathbf x^0-\lambda_{(n,0)}=(t_n,0)$. By Lemma \ref{lem:Phi-lower-bound}, there exist $N_0\in\mathbb N$ such that $|\Phi(t_n,0)|\gtrsim 1/n$ whenever $n\ge N_0$. The remainder of the proof is conducted through  the following three steps. 
	
	\underline{Step 1: \textbf{Construction of the coefficients}}. For $j\ge1$, set $N_j:=N_0 2^{j^2}$ and $I_j:=\{N_j,N_j+1,\ldots,2N_j\}$. Since
	$N_{j+1}=2^{2j+1}N_j\ge 8N_j>2N_j$, the sets $I_j$ are pairwise
	disjoint. Define $P_j:=\big(\sum_{n\in I_j}|\Phi(t_n,0)|^2\big)^{1/2}$.
	Then $P_j^2\gtrsim\sum_{n=N_j}^{2N_j}\frac1{n^2}
	\gtrsim \frac{1}{4N_j}$, and therefore $P_j\gtrsim N_j^{-1/2}$.
	We now define a sequence of coefficients
	$\alpha=(\alpha_{\mathbf k})_{\mathbf k\in\mathbb Z^2}$. For
	$n\in I_j$, let $\alpha_{(n,0)}:=\frac1j\frac{\overline{\Phi(t_n,0)}}{P_j}$,
	and set $\alpha_{\mathbf k}:=0$ for all remaining indices $\mathbf k\in\mathbb Z^2$.
	Since the sets $I_j$ are pairwise disjoint, we obtain
	\[
	\sum_{\mathbf k\in\mathbb Z^2}|\alpha_{\mathbf k}|^2=\sum_{j=1}^\infty\sum_{n\in I_j}\left|\frac1j\frac{\overline{\Phi(t_n,0)}}{P_j}\right|^2=\sum_{j=1}^\infty\frac1{j^2}\frac1{P_j^2}\sum_{n\in I_j}|\Phi(t_n,0)|^2=\sum_{j=1}^\infty \frac1{j^2}<\infty.
	\]
	Thus $\alpha\in \ell^2(\mathbb Z^2)$. By Theorem \ref{th:Riesz-bounds},  the series $f:=\sum_{\mathbf k\in\mathbb Z^2}
	\alpha_{\mathbf k}\Phi(\cdot-\lambda_{\mathbf k})$ converges in $PW_\Omega$. In particular, $f\in PW_\Omega$.
	Moreover, by the uniqueness of the Riesz-basis expansion in the Theorem \ref{thm:sampling-reconstruction}, we have $[V_{\theta,\mu}f](\lambda_{\mathbf k})=\alpha_{\mathbf k}, \mathbf k\in\mathbb Z^2$.
	
	\underline{Step 2: \textbf{A lower bound along a subsequence}}.
	Define $R_j:=h_1(N_j-1)$. Then $R_j\to\infty$ and
	$R_j\asymp N_j$. Moreover, for every $n\in I_j$,
	\[
	\|x^0-\lambda_{(n,0)}\|_2
	=
	|t_n|
	=
	h_1\left(n-\frac12\right)
	>
	h_1(N_j-1)
	=
	R_j .
	\]
	Since the only nonzero coefficients of $f$ are of the form
	$\alpha_{(n,0)}$, and
	$\alpha_{(n,0)}\Phi(t_n,0)=j^{-1}|\Phi(t_n,0)|^2/P_j\ge0$ for
	$n\in I_j$, we obtain
	\[
	\|E_{R_j}f\|_{L^\infty(\mathbb R^2)}
	\ge
	|E_{R_j}f(x^0)|
	\ge
	\sum_{n\in I_j}\alpha_{(n,0)}\Phi(t_n,0)
	=
	\frac1j\sum_{n\in I_j}\frac{|\Phi(t_n,0)|^2}{P_j}
	=
	\frac{P_j}{j}
	\gtrsim
	\frac1j N_j^{-1/2}.
	\]
	Since $R_j\asymp N_j$, this gives $\|E_{R_j}f\|_{L^\infty(\mathbb R^2)}\gtrsim
	\frac1j R_j^{-1/2}$.
	
	\underline{Step 3: \textbf{Contradiction to the faster decay}}.
	Suppose, to the contrary, that there exists $\delta>0$ such that
	$\|E_{R_j}f\|_{L^\infty(\mathbb R^2)}\le C_\delta R_j^{-1/2-\delta}$. Combining this with the lower bound obtained above gives $j^{-1}R_j^{-1/2}\lesssim R_j^{-1/2-\delta}$. Equivalently,
	$R_j^\delta/j\lesssim1$. This is impossible, since
	$R_j\asymp N_j=N_0 2^{j^2}$ and hence $R_j^\delta/j\to\infty$ as
	$j\to\infty$. Consequently, for every $\delta>0$,
	$\|E_Rf\|_{L^\infty(\mathbb R^2)}\neq O(R^{-1/2-\delta})$ as $R\to\infty$.
	This completes the proof.

	\begingroup
	\makeatletter
	\renewcommand{\@biblabel}[1]{#1.}
	\makeatother
	
	\par\medskip
	\noindent\begin{minipage}{\textwidth}
		\small
		\noindent\textsf{\textbf{Acknowledgements}}\enspace The authors thank Prof. Deguang Han at University of Central Florida for his comments on Sobolev smoothness and operator properties, which have helped improve the presentation and analysis.\par
		\medskip
		\noindent\textsf{\textbf{Funding}}\enspace The work  is  partially supported by National Natural Science Foundation of China  (No: 12471091),  Natural Science Foundation of Guangxi (No: 2026GXNSFDA00640002), and Guangxi University Outstanding Doctoral Candidate Cultivation Program Funded Project (No: ZX02010031126004).\par
	\end{minipage}
	\par
	\endgroup
\end{document}